\documentclass[a4paper]{amsart}

\usepackage{latexsym}
\usepackage{a4wide}
\usepackage{amscd}
\usepackage{graphicx}
\usepackage{amsmath}
\usepackage{amssymb}
\usepackage{mathrsfs}
\usepackage{enumerate}
\usepackage{hyperref}
\usepackage{pgf,tikz}
\usepackage[utf8]{inputenc}

\usepackage{mathtools}

\numberwithin{equation}{section} 

\newcommand{\R}{\mathbb R}
\newcommand{\N}{\mathbb N}
\newcommand{\e}{\varepsilon}

\newtheorem{theorem}{Theorem}[section]
\newtheorem{corollary}[theorem]{Corollary}
\newtheorem{lemma}[theorem]{Lemma}
\newtheorem{proposition}[theorem]{Proposition}
 \theoremstyle{definition} 
 \newtheorem{definition}[theorem]{Definition}

 \newtheorem{remark}[theorem]{Remark}

\title[Splitting of fractional eigenvalues in perforated domains]{On the
  splitting of fractional Dirichlet eigenvalues\\ in perforated
  domains}

\author[Veronica Felli]{Veronica Felli } \address{Veronica Felli,
  Dipartimento di Matematica e Applicazioni, Universit\`a degli Studi
  di Milano–Bicocca, Via Cozzi 55, 20125 Milano, Italy }
\email{veronica.felli@unimib.it}

\author[Prasun Roychowdhury]{Prasun Roychowdhury} \address{Prasun
  Roychowdhury, Department of Mathematics, Indian Institute of
  Technology Hyderabad, Kandi, Sangareddy, Telangana, 502285, India }
\email{prasunrc@math.iith.ac.in}

\author[Giovanni Siclari]{Giovanni Siclari} \address{Giovanni Siclari,
  Centro di Ricerca Matematica Ennio De Giorgi, Scuola Normale
  Superiore di Pisa, Piazza dei Cavalieri 3, 56126 Pisa, Italy }
\email{giovanni.siclari@sns.it}

\date{November 27, 2026}

\begin{document}

\begin{abstract}
  In this paper, we study the asymptotic behavior of eigenvalues of
  the fractional Laplacian under singular domain perturbations
  produced by the removal of a hole with vanishing fractional
  capacity. First, we establish an expansion formula for eigenbranches
  bifurcating from a possibly multiple eigenvalue.  We then
  investigate splitting phenomena: specifically, removing a family of
  holes concentrating at a zero-capacity compact set (lying outside 
  a closed set of measure zero) causes a multiple eigenbranch to split into
  branches of smaller multiplicity. Finally, for holes concentrating
  at a point, we prove that the largest eigenbranch generically
  detaches from the rest and becomes simple.
\end{abstract}

\maketitle


\medskip 


\section{Introduction}
 This paper is concerned with describing how the eigenvalues of the
fractional Laplacian, in a bounded open set with Dirichlet boundary
conditions, respond to the removal of small holes. In particular, we
aim to study any splitting phenomenon of multiple eigenvalues produced
by this type of perturbation.

In spectral theory and partial differential equations, the sensitivity
of the spectrum of a differential operator to small perturbations
attracts considerable interest, driving a rich line of
research. Within this topic, questions regarding the genericity of
simple eigenvalues versus the persistence of multiplicity have been of
particular importance ever since the pioneering contributions in the
1970s by Micheletti
\cite{micheletti1972-2,micheletti1972-1,micheletti1973,micheletti1976},
Uhlenbeck \cite{uhlenbeck1972,uhlenbeck1976}, and Albert
\cite{albert1975,albert1978}.  The research presented in the
aforementioned papers supports the conclusion that the simplicity of
the spectrum of a differential operator is a generic property, in the
sense that, for an operator with multiple eigenvalues, almost any
sufficiently small perturbation produces an operator with only simple
eigenvalues. It is worth noting that these works deal with regular
perturbations, meaning that the perturbed problem is obtained by
applying a diffeomorphism to the original operator, with perturbations
potentially affecting both the coefficients of the operator (see
\cite{albert1975,albert1978,uhlenbeck1972,uhlenbeck1976}) and the
domain (see
\cite{micheletti1972-2,micheletti1972-1,micheletti1973,micheletti1976}).

In the present paper, we are instead concerned with singular
perturbations, i.e., those not arising from the application of regular
transformations. A typical example (and indeed the one on which we
will focus) is given by the perturbation caused by the removal of a
small compact hole $K$ from the interior of a bounded open set
$\Omega$: this is a singular perturbation, as $\Omega$ and
$\Omega\setminus K$ are not diffeomorphic. This type of perturbation
has been widely studied in the literature for the Laplacian,
particularly under Dirichlet boundary conditions. The appropriate
notion of smallness for the removed hole that guarantees spectral
stability for the Dirichlet Laplacian is related to the classical
notion of capacity, as observed in \cite{rauch-taylor}. In particular,
removing a hole of null capacity from the domain leaves the spectrum
unchanged. Consequently, a perturbation theory can be developed using
the capacity of the removed set as a perturbation parameter, showing
that every eigenvalue of the problem in the perforated domain
converges to the corresponding unperturbed eigenvalue, as the capacity
of the hole goes to zero; see \cite{courtois} for manifolds without
boundary and \cite{Bertrand-Colbois} for manifolds with
boundary. Furthermore, in \cite{courtois}, Courtois obtains asymptotic
expansions for the eigenvalue variation in terms of $u$-capacities,
i.e. of a modified notion of capacity incorporating limit
eigenfunctions $u$ on the hole.  More precisely, for a simple
eigenvalue, the variation is asymptotic to the $u$-capacity relative
to a suitably normalized eigenfunction $u$ (see also \cite{AFHL}),
whereas, in the multiple case, one obtains an expansion in terms of
the integral of the squared eigenfunctions with respect to a suitable
equilibrium measure associated with the hole (see also
\cite{ALM}). For the Dirichlet Laplacian, of particular interest is
also the paper \cite{flucher}, which provides asymptotic expansions
for the eigenvalues in a domain from which
a hole concentrating at a point is removed, observing a splitting
phenomenon (see \cite[Theorem 6]{flucher}). 

For the Neumann Laplacian, eigenvalue stability in the presence of
small holes requires different conditions compared to the Dirichlet
case. In particular, by \cite{rauch-taylor} this is guaranteed under
the assumption that there exist extension operators for $H^1$ Sobolev
spaces with uniformly bounded operator norms; this property holds,
e.g., for holes given by a fixed shape rescaling to a point. For the
Neumann case, in addition to \cite{rauch-taylor}, we mention
\cite{jimbo,LdC2012, LNS2011,ozawa1983} and the more recent
contributions \cite{FLO_1} for simple eigenvalues and \cite{FLO_2} for
multiple ones, showing that a kind of boundary torsional rigidity
replaces capacity as the perturbation parameter in the Neumann
setting; see also \cite{FRS} for Robin boundary
conditions. We also mention \cite{BORS}, where a
  unified abstract framework in a Hilbertian setting is developed to
  obtain quantitative asymptotic expansions for eigenvalue variations
  under singular perturbations, covering both simple and multiple
  eigenvalues.

In the present paper, we aim to investigate splitting phenomena for
multiple Dirichlet eigenvalues of the fractional Laplacian
$(-\Delta)^s$ arising from the singular perturbation induced by the
removal of a small hole. Regarding regular perturbations in the
fractional setting, we mention that generic properties for the
eigenvalues of the fractional Dirichlet Laplacian were studied in
\cite{FGMP}.  Employing an abstract result by Micheletti
\cite{micheletti1976}, \cite{FGMP} shows in particular that, given a
smooth bounded domain $\Omega$, there exists an arbitrarily small
regular perturbation of $\Omega$ where all fractional Dirichlet
eigenvalues are simple; moreover, in \cite{FGMP}, perturbations of the
coefficients are also considered.  On the other hand, the persistence
of the multiplicity of eigenvalues of the fractional Laplacian is
studied in \cite{GMP}, under regular perturbations of either the
domain or the coefficients of the operator.  Related results are
contained in \cite{DFW}, where a Hadamard formula for the derivative
of the Dirichlet eigenvalues of $(-\Delta)^s$, with respect to domain
deformations, is provided and used to prove the simplicity of radial
eigenvalues in a ball or an annulus.

Regarding the effect of singular perturbations, specifically the
removal of small holes from the domain, on the eigenvalues of the
fractional Laplacian, a key reference is \cite{AFN}.  In particular,
in \cite{AFN}, a notion of fractional capacity is introduced (see
Definition \ref{def:s-capacity}), spectral stability is proved as the
fractional capacity of the hole goes to zero, and asymptotic
expansions of the variation of simple eigenvalues are
obtained. Capacitary potentials and equilibrium measures for
fractional energies, as well as the relation between zero capacity and
removability of singular sets for nonlocal harmonic functions, are
also investigated in \cite{FocPalZep26,Pal26a}.

The purpose of this manuscript is twofold. On the one hand, we prove
an asymptotic expansion formula for the eigenbranches of the
fractional Laplacian that bifurcate from a possibly multiple
eigenvalue, as a hole with vanishing fractional capacity is removed
from the domain. This is achieved in our first main result, Theorem
\ref{theo_main}, which generalizes the expansion obtained in
\cite[Theorem 1.5]{AFN} for simple eigenvalues to the multiple case,
as well as providing the fractional counterpart of Courtois' result
\cite[Theorem 1.2]{courtois} for the classical Laplacian. On the other
hand, we address genericity issues, using the expansion proved in
Theorem \ref{theo_main} to investigate splitting phenomena for
multiple eigenbranches. First, in Theorem \ref{t:splitting}, we show
that, by removing a family of holes concentrating at a compact set of
zero capacity, which lies outside a closed set of measure zero, the
multiple eigenbranch splits into branches of smaller
multiplicity. Then, in Theorem \ref{cor_main}, we consider the special
case of a family of holes concentrating at a single point; if the
point lies outside the nodal set of at least one eigenfunction (which
occurs almost everywhere), it is precisely the largest eigenbranch
that detaches from all the others, becoming simple. In particular,
Theorem \ref{cor_main} provides the fractional counterpart to
Flucher's result \cite[Theorem 6]{flucher}.

In the next section, after formulating the problem and introducing the
necessary notation, we state the main results of the paper.

\section{Statement of the main results}
In this section, we first set up the problem and introduce the
required definitions and notation, before presenting our main results.

\subsection{Statement of the problem}
Let $N\geq1$ and $s\in(0,\min\{1,N/2\})$. For any open set
$\Omega\subset\R^N$, we consider the homogeneous fractional Sobolev
space $\mathcal D^{s,2}(\Omega)$ defined as the completion of
$C^\infty_c(\Omega)$ with respect to the Gagliardo norm
\[
[u]_{H^s(\R^N)}=\left(\int_{\R^N}\int_{\R^N}
  \frac{|u(x)-u(y)|^2}{|x-y|^{N+2s}} \,dx\,dy\right)^{1/2},
\]
where, in the double integral above, the function $u$ is understood to
be trivially extended to $\R^N\setminus\Omega$.  We equip
$\mathcal D^{s,2}(\Omega)$ with the following scalar product (which makes it a Hilbert space)
\begin{align}\label{eq:21}
  (u,v)_{\mathcal
  D^{s,2}(\Omega)}&=\frac{C(N,s)}2\int_{\R^{N}}\int_{\R^{N}}
                    \frac{(u(x)-u(y))(v(x)-v(y))}{|x-y|^{N+2s}}\,dx\,dy\\
  &\notag
  =\int_{\R^N}
  |\xi|^{2s}\mathcal F u(\xi)\overline{\mathcal F v(\xi) }\,d\xi,
\end{align}
and the associated norm 
\[
\|u\|_{\mathcal D^{s,2}(\Omega)}=(u,u)_{\mathcal
    D^{s,2}(\Omega)}^{1/2}=\sqrt{\tfrac{C(N,s)}2}\,[u]_{H^s(\R^N)},
\]
where 
\begin{equation}\label{eq:14}
C(N,s)=\pi^{-\frac
  N2}2^{2s}\frac{\Gamma\big(\frac{N+2s}{2}\big)}{\Gamma(2-s)}s(1-s),
\end{equation}
$\Gamma$ is the Gamma function, and 
$\mathcal F u$ denotes the unitary Fourier transform of $u$.

The following notions of \emph{fractional capacities} play a pivotal
role in the perturbation theory developed in \cite{AFN} for simple
eigenvalues of singularly perturbed fractional problems. For an
overview of properties of the Gagliardo $s$-fractional capacity, we
refer the reader to \cite[Appendix A]{ritorto}; related concepts of
fractional capacity are explored in
\cite{adams,Adams-Xiao,shi-xiao,Warma15}.

\begin{definition}\label{def:s-capacity}
Let $\Omega\subset\R^N$ be a bounded open set and  $s\in
(0,\min\{1,N/2\})$. Let $K\subset\Omega$ be
a compact set and let $\zeta_K\in C^\infty_c(\Omega)$ be such that
$\zeta_K\equiv1$ in a neighborhood of $K$.
\begin{enumerate}[(i)]
\item The \emph{Gagliardo $s$-fractional capacity} of $K$ in $\Omega$
  is defined  as
\begin{equation}\label{eq:def-capOK}
\mathop{\rm Cap}\nolimits^s_\Omega(K)=\inf \left\{
\|w\|_{\mathcal D^{s,2}(\Omega)}^2: \, w\in \mathcal D^{s,2}(\Omega)  \text{ and }  
w-\zeta_K\in \mathcal D^{s,2}(\Omega\setminus K) \right\}.
\end{equation}
\item
For every $u\in
  H^s_{\rm loc}(\Omega)$,  the \emph{$s$-fractional
    $u$-capacity} of $K$ in $\Omega$ is  defined  as
\begin{equation*}
\mathop{\rm Cap}\nolimits^s_\Omega(K,u)=
\inf \left\{
\|w\|_{\mathcal D^{s,2}(\Omega)}^2: \, w\in \mathcal D^{s,2}(\Omega)  \text{ and }  
w-\zeta_K u\in \mathcal D^{s,2}(\Omega\setminus K) \right\}.
\end{equation*}
If  $u\in \mathcal D^{s,2}(\Omega)$, the $s$-fractional
    $u$-capacity of $K$ in $\Omega$ can also be  characterized as
\begin{equation*}
\mathop{\rm Cap}\nolimits^s_\Omega(K,u)=
\inf \left\{
\|w\|_{\mathcal D^{s,2}(\Omega)}^2: \, w\in \mathcal D^{s,2}(\Omega)  \text{ and }  
w-u\in \mathcal D^{s,2}(\Omega\setminus K) \right\}.
\end{equation*}
\end{enumerate}
\end{definition}
We are interested in the following eigenvalue problem for the
\emph{restricted fractional Laplacian} in an open bounded set
$\Omega\subset\R^N$, under homogeneous Dirichlet boundary conditions,
\begin{equation}\label{eq:eigenvalue_main}
\begin{cases}
(-\Delta)^s u=\lambda u, \quad& \text{in } \Omega, \\
u=0, &\text{in }\R^N\setminus\Omega.
\end{cases}
\end{equation}
We recall that the operator $(-\Delta)^s$ is defined over the space
$C^\infty_c(\R^N)$ by the principal value integral
\begin{equation*}
(-\Delta)^s u(x)=C(N,s) \lim_{\e\to0^+} \int_{\R^N\setminus B_\e(x)} 
\frac{u(x)-u(y)}{|x-y|^{N+2s}} \,dy,
\end{equation*}
where $C(N,s)$ is given in \eqref{eq:14}, or equivalently through the
Fourier transform:
\begin{equation*}
\mathcal{F}((-\Delta)^s u)(\xi)=|\xi|^{2s} \mathcal{F}u(\xi), \qquad \xi\in \R^N.
\end{equation*}
The operator $(-\Delta)^s$ extends to a bounded linear map from 
$\mathcal D^{s,2}(\R^N)$ into its dual space $(\mathcal{D}^{s,2}(\R^N)
)^*$. In fact, this extension coincides with the Riesz isomorphism for
$\mathcal{D}^{s,2}(\R^N)$ induced by the scalar product
\eqref{eq:21}, namely
\begin{equation*}
  \phantom{a}_{( \mathcal{D}^{s,2}(\R^N) )^*}\langle(-\Delta)^s u,
  v\rangle_{ \mathcal{D}^{s,2}(\R^N) } = (u,v)_{\mathcal{D}^{s,2}(\R^N)}
\end{equation*}
for every $u,v\in \mathcal{D}^{s,2}(\R^N)$.
Hence, the eigenvalue problem \eqref{eq:eigenvalue_main}  can be
understood in the following weak sense: $\lambda\in\R$ is an eigenvalue of
\eqref{eq:eigenvalue_main} 
if there exists some eigenfunction $u\in \mathcal
  D^{s,2}(\Omega)\setminus\{0\}$ such that
  \begin{equation}\label{eq_eigen}
 (u,v)_{\mathcal
    D^{s,2}(\Omega)}=\lambda \int_{\Omega} u(x)v(x) \,dx \quad
\text{for every }  
v\in \mathcal D^{s,2}(\Omega).
\end{equation}
By classical spectral theory, \eqref{eq:eigenvalue_main} admits a
diverging sequence
of eigenvalues
 with finite multiplicity 
 \begin{equation*}
0< \lambda_1^s(\Omega)\leq \lambda_2^s(\Omega) \leq
\lambda_3^s(\Omega)
\leq \ldots \to+\infty.
\end{equation*}
In \cite{AFN} it was proved that, if $\Omega\subset\R^N$ is a bounded
open set and $K\subset\Omega$ is compact, then, for every
$k\in\N\setminus\{0\}$,
\begin{equation}\label{limit_eigen}
\lambda_k^s(\Omega\setminus
K)\to\lambda_k^s(\Omega)\quad\text{as }\mathop{\rm
  Cap}\nolimits^s_\Omega(K)\to 0,
\end{equation}
i.e., the spectrum is stable under perturbations arising from the
excision of a small hole, provided its capacity approaches zero.
Furthermore, \cite{AFN} shows that the variation of a simple
eigenvalue due to perforation is asymptotic to the
$s$-fractional $u$-capacity of the removed set, where $u$ is a
normalized eigenfunction associated with the unperturbed eigenvalue,
as the hole concentrates at a set of zero $s$-fractional capacity in
the sense of Definition \ref{def:concentrates}.

From an energy perspective, the fact that removing a set of small
capacity only slightly perturbs the spectrum (and, in the limit,
removing sets of zero capacity leaves the spectrum unchanged) is
closely related to the Serrin-type removability result obtained in
\cite[Theorem 1.1]{Pal26a}. There, it is shown (in the more general
setting of $(s,p)$-capacities) that a set is removable for harmonic
functions if and only if it has zero capacity; see also \cite{kimlee}
for the sufficiency of this condition.

The present paper aims to quantify the spectral stability stated in
\eqref{limit_eigen}, including the case of possibly multiple
eigenvalues, as $\mathop{\rm Cap}\nolimits^s_\Omega(K)\to 0$, with the
specific goal of observing and investigating potential multiple-branch
bifurcation phenomena.

Before stating our main results, we briefly
review additional properties of the fractional capacity
$\mathop{\rm Cap}\nolimits^s_\Omega(K)$.

\subsection{Capacitary potentials and equilibrium measures}\label{subsec_equil_measu}

By a standard variational argument, for any compact set
$K \subset \Omega$, the infimum in \eqref{eq:def-capOK} is
  achieved by a unique minimizer $v_{\Omega,K}$, which also turns out to
  be the unique weak solution to the problem
\begin{equation*}
\begin{cases}
(-\Delta)^s v_{\Omega,K}=0, \quad& \text{in } \Omega\setminus K, \\
v_{\Omega,K}=0, &\text{in }\R^N\setminus\Omega,\\
v_{\Omega,K}=1, &\text{in }K,
\end{cases}
\end{equation*}
in the sense that, letting $\zeta_K$ be as in Definition
  \ref{def:s-capacity},
$v_{\Omega,K}-\zeta_K\in \mathcal D^{s,2}(\Omega\setminus K)$
and
\begin{equation}\label{eq:eqvok}
  (v_{\Omega,K},\phi)_{\mathcal
    D^{s,2}(\R^N)}=0\quad\text{for every }\phi \in \mathcal
  D^{s,2}(\Omega\setminus K).
\end{equation}
In particular,
\begin{equation*}
\mathop{\rm Cap}\nolimits^s_\Omega(K)=\|v_{\Omega,K}\|^2_{\mathcal{D}^{s,2}(\Omega)}.
\end{equation*}
We will refer to $v_{\Omega,K}$ as the \emph{capacitary potential} of $K$ in
$\Omega$.  

We observe that, if $\phi\in \mathcal{D}^{s,2}(\Omega)$, then
    its truncation
    \begin{equation*}
      T\phi:=\min\{\phi,1\}=1-(1-\phi)^+
    \end{equation*}
    satisfies
    \begin{equation*}
      T\phi\in \mathcal{D}^{s,2}(\Omega)\quad
      \text{and}\quad 
      \|T\phi\|_{\mathcal{D}^{s,2}(\Omega)}\leq
      \|\phi\|_{\mathcal{D}^{s,2}(\Omega)}.
    \end{equation*}
    Furthermore, if $\phi-\zeta_K\in \mathcal{D}^{s,2}(\Omega\setminus
    K)$, then
    $T\phi-\zeta_K\in \mathcal{D}^{s,2}(\Omega\setminus
    K)$. Therefore, the capacity $\mathop{\rm
      Cap}\nolimits^s_\Omega(K) $ in \eqref{eq:def-capOK} can be
    equivalently characterized as
\begin{equation*}
\mathop{\rm Cap}\nolimits^s_\Omega(K)=\inf \left\{
\|w\|_{\mathcal D^{s,2}(\Omega)}^2: \, w\in \mathcal D^{s,2}(\Omega)  \text{ and }  
Tw-\zeta_K\in \mathcal D^{s,2}(\Omega\setminus K) \right\},
\end{equation*}
see also \cite[Section 2]{Pal26a} and \cite[Proposition A.5]{FocPalZep26}.
It follows that, for every $\varphi\in C^\infty_c(\Omega)$ such that
$\varphi\geq 0$ in $\Omega$,
$\|v_{\Omega,K}\|^2_{\mathcal{D}^{s,2}(\Omega)}\leq
\|v_{\Omega,K}+\varphi\|^2_{\mathcal{D}^{s,2}(\Omega)}$, so that
$(v_{\Omega,K},\varphi)_{\mathcal{D}^{s,2}(\Omega)}\geq0$.
In particular, the map
\begin{equation*}
  \Lambda_{\Omega,K}: C^\infty_c(\Omega)\to \R,\quad
  \Lambda_{\Omega,K}(\varphi)= (v_{\Omega,K},\varphi)_{\mathcal{D}^{s,2}(\Omega)}
\end{equation*}
defines a distribution which is non-negative on non-negative test
functions, i.e. $\Lambda_{\Omega,K}(\varphi)\geq0$ for every
$\varphi\in C^\infty_c(\Omega)$ such that $\varphi\geq 0$ in $\Omega$.
By the Riesz-Schwartz Theorem, see \cite[Theorem
10.13]{L_book_sobolev}, there exists a unique
Radon measure $\mu_{\Omega,K}^s: \mathcal B(\Omega) \to[0,+\infty]$ such that
\begin{equation}\label{eq:Lmu}
 \Lambda_{\Omega,K}(\varphi)=
 (v_{\Omega,K},\varphi)_{\mathcal{D}^{s,2}(\Omega)}=\int_\Omega\varphi\,d\mu_{\Omega,K}^s
 \quad\text{for
   all }\varphi\in C^\infty_c(\Omega).
\end{equation}
We observe that the support of $\mu_{\Omega,K}^s$ is contained in
$K$.

Furthermore, since
$v_{\Omega,K}-\zeta_K\in \mathcal D^{s,2}(\Omega\setminus K)$, there
exists a sequence $\varphi_n\in C^\infty_c(\Omega\setminus K)$ such
that $\varphi_n \to v_{\Omega,K}-\zeta_K$ in
$\mathcal D^{s,2}(\R^N)$. Hence, by \eqref{eq:eqvok},
\begin{align}\label{eq_muK}
\mathop{\rm Cap}\nolimits^s_\Omega(K)&=(v_{\Omega,K} ,v_{\Omega,K}
)_{\mathcal{D}^{s,2}(\R^N)}=
\lim_{n\to\infty}(v_{\Omega,K} ,\varphi_n+\zeta_K
                                       )_{\mathcal{D}^{s,2}(\R^N)}\\
  &\notag=(v_{\Omega,K} ,\zeta_K
)_{\mathcal{D}^{s,2}(\R^N)}=\int_\Omega \zeta_K \,d\mu_{\Omega,K}^s=
\int_K \zeta_K \,d\mu_{\Omega,K}^s=\mu_{\Omega,K}^s(K).
\end{align}
We also observe that, since the support of $\mu_{\Omega,K}^s$ is contained in
$K$, we have
\begin{equation*}
  \mu_{\Omega,K}^s(\Omega)= \mu_{\Omega,K}^s(K)= \mathop{\rm Cap}
  \nolimits^s_\Omega(K).
\end{equation*}
We refer to $\mu_{\Omega,K}^s$ as the \textit{$s$-equilibrium measure}
of $K$ in $\Omega$.  We point out that $\mu_{\Omega,K}^s$ is a
particular case of the family of equilibrium measures associated with
relative $(s,p)$-capacities, with symmetric measurable kernels
comparable to the Gagliardo one, constructed in \cite[Lemma
3.1]{Pal26a}; see also \cite[Lemma 5.1 and Proposition
6.2]{FocPalZep26} for the counterpart on the whole space.

\subsection{Main results}

Our first result is an asymptotic expansion of
$\lambda_{k}^{s}(\Omega\setminus K_n)$ for a sequence of compact sets
$\{K_n\}_{n\in\N}$ contained in $\Omega$, whose $s$-fractional
capacity tends to zero, which remain well-separated from the boundary
$\partial\Omega$, and whose normalized $s$-equilibrium measures converge
weakly in the sense of measures.

  In particular, we assume that
  \begin{align*}
    &\label{eq:K1}\tag{K1}   \text{for every $n\in\N$, $K_n$ is a compact set  such that $K_n\subset\Omega$
      and $\mathop{\rm Cap}\nolimits^s_\Omega(K_n)>0$},\\
    &\label{eq:K2}\tag{K2}   \lim_{n\to\infty}
      \mathop{\rm Cap}\nolimits^s_\Omega(K_n)=0,\\
    &\label{eq:K3}\tag{K3}\text{there exists $\delta>0$ such that
      $\mathop{\rm dist}(K_n,\partial\Omega)>\delta$ for every $n\in\N$},\\
    &\label{hp_weak_conv_mu}\tag{K4}\text{there exists a 
      Radon measure
      $\tilde \mu: \mathcal B(\Omega) \to[0,+\infty]$ such that}
      \shortintertext{\centering 
      $\tilde{\mu}^s_{\Omega,K_n}:=\dfrac{\mu_{\Omega,K_n}^s}
      {\mu_{\Omega,K_n}^s(\Omega)}=
      \dfrac{\mu_{\Omega,K_n}^s}
      {\mathop{\rm Cap}
      \nolimits^s_\Omega(K_n)}\rightharpoonup \tilde
      \mu$ as $n\to\infty$,}\\[-10pt]
    \notag  &\text{where $\mu_{\Omega,K_n}^s$ is the
              $s$-equilibrium measure of $K_n$ in $\Omega$ introduced in Section
              \ref{subsec_equil_measu}.} 
    \end{align*}
    For our purposes, the assumption in \eqref{eq:K1} that
    $\mathop{\rm Cap}\nolimits^s_\Omega(K_n)>0$ for every $n\in\N$ is
    not restrictive; indeed, if
    $\mathop{\rm Cap}\nolimits^s_\Omega(K_n)=0$ for some $n$, we have
    $\lambda_k^s(\Omega\setminus K_n)=\lambda_k^s(\Omega)$ for all
    $k$, making the expansion we intend to prove trivial. Furthermore,
    by \cite{AFN}, assumption \eqref{eq:K2} ensures spectral stability, i.e.
    $\lim_{n\to\infty}\lambda_k^s(\Omega\setminus
    K_n)=\lambda_k^s(\Omega)$ for all $k\geq1$.

    The convergence in \eqref{hp_weak_conv_mu} is understood in the
    weak sense of measures, meaning that
\begin{equation*}
  \lim_{n\to\infty}\int_{\Omega} \phi \, d\tilde{\mu}^s_{\Omega,K_n}
  =\int_{\Omega} \phi \, d\tilde{\mu} \quad \text{for any } \phi
  \in C^0_c(\Omega).
\end{equation*}
Assumption \eqref{eq:K3} implies that the sequence
$\{\tilde{\mu}^s_{\Omega,K_n}\}_{n\in\N}$ is tight, thus ensuring that
the weak limit $\tilde\mu$ in \eqref{hp_weak_conv_mu} is necessarily a
probability measure.

Let $n_0,m_0 \in \N\setminus\{0\}$ be such that
$\lambda^s_{n_0}(\Omega)$ has multiplicity $m_0$ as an eigenvalue of
\eqref{eq_eigen}, with
\begin{equation}\label{lambda_mul_m}
  \lambda^s_{n_0-1}(\Omega)< \lambda^s_{n_0}(\Omega)=
  \lambda^s_{n_0+1}(\Omega)=\dots= \lambda^s_{n_0+m_0-1}(\Omega)
  <\lambda^s_{n_0+m_0}(\Omega).
\end{equation}
Let $E(\lambda^s_{n_0}(\Omega))\subset \mathcal D^{s,2}(\Omega)$ be
the corresponding eigenspace.  Let us define the symmetric bilinear
form
\begin{equation}\label{eq:bilform}
r_{\tilde\mu}:E(\lambda^s_{n_0}(\Omega))\times E(\lambda^s_{n_0}(\Omega)) \to
\R,  \quad  r_{\tilde\mu}(u,v):=\int_\Omega u v \, d\tilde \mu \quad 
\text{ for any } u,v \in E(\lambda^s_{n_0}(\Omega)),
\end{equation}
and denote its eigenvalues by $\{\beta_i\}_{i=1,\dots,m_0}$ (in
increasing order) counted with multiplicity.

\begin{theorem}\label{theo_main}
  Let $N\geq1$, $s\in(0,\min\{1,N/2\})$, and $\Omega\subset\R^N$ be a
  bounded open set.  Let $n_0,m_0 \in \N\setminus\{0\}$ be such that
  $\lambda^s_{n_0}(\Omega)$ has multiplicity $m_0$ as an eigenvalue of
  \eqref{eq_eigen}, with the index $n_0$ chosen as in
  \eqref{lambda_mul_m}.  Let $\{K_n\}_{n\in\N}$ be a sequence of
  compact sets contained in $\Omega$ satisfying \eqref{eq:K1},
  \eqref{eq:K2}, \eqref{eq:K3}, and \eqref{hp_weak_conv_mu}.  Then,
  for every $1\leq i\leq m_0$,
\begin{equation*}
  \lambda_{n_0+i-1}^{s}(\Omega\setminus K_n) =\lambda^s_{n_0} (\Omega)+
  \beta_i\mathop{\rm Cap}\nolimits^s_\Omega(K_n)+
  o\left(\mathop{\rm Cap}\nolimits^s_\Omega(K_n)\right)\quad\text{as } n\to\infty,
\end{equation*}
where $\{\beta_i\}_{i=1,\dots,m_0}$ are the
eigenvalues (in increasing order) of the bilinear form $r_{\tilde\mu}$
defined in \eqref{eq:bilform}.
\end{theorem}

By Theorem \ref{theo_main}, if $m_0>1$ and there exists
$j \in \{2,\dots,m_0\}$ such that $\beta_{j-1}< \beta_j$, then a
splitting phenomenon occurs as a result of the domain perforation,
i.e. the multiple eigenvalue $\lambda^s_{n_0} (\Omega) $ splits into
more (at least two) distinct branches.  To identify situations in
which such a bifurcation occurs, we consider the special case of a
family of compact sets concentrating at a limit compact set of zero
capacity.

\begin{definition}\label{def:concentrates}
  Let $\Omega\subset\R^N$ be a bounded open set.  Let
  $\{K(\varepsilon)\}_{\varepsilon>0}$ be a family of compact sets
  contained in $\Omega$.  We say that $K(\varepsilon)$ is
  concentrating at a compact set $K\subset\Omega$ as $\e\to0$, if, for
  every open set $\omega$ such that $K\subset\omega\subseteq \Omega$,
  there exists $\varepsilon_\omega>0$ such that
  $K(\varepsilon)\subset\omega$ for every
  $0<\varepsilon<\varepsilon_\omega$.
\end{definition}
We observe that, if $\{K(\varepsilon)\}_{\varepsilon>0}$ is concentrating
  at a compact set $K\subset\Omega$  with $\mathop{\rm
    Cap}\nolimits^s_\Omega(K)=0$, then, by \cite[Lemma 4.2]{AFN},
  \begin{equation*}
    \lim_{\e\to 0}\mathop{\rm
    Cap}\nolimits^s_\Omega(K(\e))=0.
\end{equation*}
Moreover, there exists some compact set $\bar K\subset\Omega$ such
that $K(\e)\subset\bar K$ for $\e$ sufficiently small. On the other
hand, by \cite[Theorem 1.41]{EG_book}, for any sequence $\e_n\to0^+$,
there exists a subsequence $\{\e_{n_k}\}$ such that
$\tilde{\mu}^s_{\Omega,K(\e_{n_k})}$ weakly converges to some
probability Radon measure $\tilde\mu$ with support contained in $K$.
Note that such a limiting measure $\tilde\mu$, being supported on a
set of zero capacity, can only emerge as a genuinely singular weak
limit; indeed, as observed in \cite[Lemma 6.1]{FocPalZep26},
equilibrium measures do not charge polar sets.

Theorem \ref{theo_main} then applies to the sequences
$\{K(\e_{n_k})\}_k$, providing an expansion of the eigenvalues
$\lambda_{n_0+i-1}^{s}(\Omega\setminus K(\e_{n_k}))$,
$i=1,\dots,m_0$. By studying when the coefficients $\beta_i$ appearing
in the expansion do not overlap for different indices, we can
establish the following splitting result for compact concentrating
sets.

\begin{theorem}\label{t:splitting}
  Let $N\geq1$, $s\in(0,\min\{1,N/2\})$, and $\Omega\subset\R^N$ be a
  bounded open set.  Let $n_0,m_0 \in \N\setminus\{0\}$ be such that
  $\lambda^s_{n_0}(\Omega)$ has multiplicity $m_0\ge 2$ as
  an eigenvalue of \eqref{eq_eigen}, with the index $n_0$ chosen as in
  \eqref{lambda_mul_m}.  Let
  $\{u_{n_0}, u_{n_0+1},\dots, u_{n_0+m_0-1}\}$ be an
  $L^2(\Omega)$-orthonormal basis of the associated $m_0$-dimensional
  eigenspace.  Let $\{K(\varepsilon)\}_{\varepsilon>0}$ be a family of
  compact sets contained in $\Omega$ with
  $\mathop{\rm Cap}\nolimits^s_\Omega(K(\e))>0$ for every $\e>0$,
  concentrating at a compact set $K\subset\Omega$ such that
  \begin{equation*}
    \mathop{\rm
    Cap}\nolimits^s_\Omega(K)=0
\end{equation*}
and
\begin{equation}\label{eq:ass-conn-comp}
  \text{$K$ is contained in a single connected component of
    $\Omega\setminus \mathcal N$},
\end{equation}
where $\mathcal N=\bigcup_{i=1}^{m_0}u_{n_0+i-1}^{-1}(0)$. Then, there
exists $\e_0>0$ such that, for every $\e\in(0,\e_0)$ and
$i=1,\dots,m_0$,
\begin{equation*}
  m\left(\lambda_{n_0+i-1}^{s}(\Omega\setminus K(\e))\right)< m_0,
\end{equation*}
where $m\left(\lambda_{n_0+i-1}^{s}(\Omega\setminus K(\e))\right)$
denotes the multiplicity of the eigenvalue
$\lambda_{n_0+i-1}^{s}(\Omega\setminus K(\e))$.
\end{theorem}
We note that, by Schauder estimates for the fractional Laplacian (see,
e.g., \cite[Theorem 2.15]{TLX}), the eigenfunctions $u_{n_0+i-1}$ are
of class $C^\infty(\Omega)$. Moreover, it is known from \cite[Theorem
1.4]{FallFelli14} that the nodal sets of the eigenfunctions have zero
Lebesgue measure. Thus, the set $\mathcal N$ appearing in Theorem
\ref{t:splitting} is a closed set in $\Omega$ of zero Lebesgue
measure, making the condition that $K(\e)$ concentrates at a set lying
in a connected component of the complement both possible and quite
generic.

The above result can be made more precise in the case where the limit
compact set reduces to a singleton $K=\{x_0\}$ for some
$x_0 \in \Omega$; in this case, $\tilde\mu=\delta_{x_0}$ and
the convergence of normalized $s$-equilibrium measures occurs not only
along subsequences, thus yielding the following theorem.

\begin{theorem}\label{cor_main}
  Let $N\geq1$, $s\in(0,\min\{1,N/2\})$, and $\Omega\subset\R^N$ be a
  bounded open set.  Let $n_0,m_0 \in \N\setminus\{0\}$ be such that
  $\lambda^s_{n_0}(\Omega)$ has multiplicity $m_0$ as an eigenvalue of
  \eqref{eq_eigen}, with the index $n_0$ chosen as in
  \eqref{lambda_mul_m}.  
  Let $\{u_{n_0}, u_{n_0+1},\dots, u_{n_0+m_0-1}\}$ be an
  $L^2(\Omega)$-orthonormal basis of the associated $m_0$-dimensional
  eigenspace.
 Let $x_0 \in \Omega$ and let
  $\{K(\varepsilon)\}_{\varepsilon>0}$ be a family of compact sets
  contained in $\Omega$ concentrating at $\{x_0\}$ in the sense of
  Definition \ref{def:concentrates}.  Then
\begin{equation}\label{eq:cor_main1}
  \lambda_{n_0+m_0-1}^{s}(\Omega\setminus K(\e)) =\lambda^s_{n_0}
  (\Omega)+
  \left(\sum_{i=1}^{m_0}u_{n_0+i-1}^2(x_0)\right)
  \mathop{\rm Cap}\nolimits^s_\Omega(K(\e))+o\left(\mathop{\rm Cap}\nolimits^s_\Omega(K(\e))\right)
\end{equation}
as $\e \to0^+$. Furthermore, if $m_0>1$,
\begin{equation}\label{eq:cor_main2}
  \lambda_{n_0+i-1}^s(\Omega\setminus K(\e)) =\lambda_{n_0}^s(\Omega)
  +o\left(\mathop{\rm Cap}\nolimits^s_\Omega(K(\e))\right)
\end{equation}
as $\e \to0^+$, for every $1\leq i\leq m_0-1$.
\end{theorem}
From Theorem \ref{cor_main}, it follows that, if $x_0 \in \Omega$ is
such that at least one eigenfunction of problem \eqref{eq_eigen}
associated with the eigenvalue $\lambda_{n_0}^s(\Omega)$ is not zero
at $x_0$, then, for $\e$ small enough, the eigenvalue
$\lambda_{n_0+m_0-1}^s(\Omega\setminus K(\e))$ is simple, while the
eigenvalues $\{\lambda_{n_0+i-1}^s(\Omega\setminus K(\e))\}_{i<m_0}$
have multiplicity strictly smaller than $m_0$. In particular, a
splitting phenomenon occurs, i.e.  perforating the domain causes the
multiple eigenvalue $\lambda^s_{n_0} (\Omega)$ to split into at least
two distinct branches. We note that the condition $u(x_0) \neq 0$ for
an eigenfunction $u$
of \eqref{eq_eigen} is \textit{generic},
in the sense that $u(x) \neq 0$ for a.e. $x \in \Omega$, as proved in
\cite[Theorem 1.4]{FallFelli14}.

The point-concentration regime considered in Theorem \ref{cor_main} is
reminiscent of \cite[Proposition 7.2]{FocPalZep26}, where, for
shrinking balls, one observes the convergence of the normalized
equilibrium measures to $\delta_{x_0}$, together with the convergence
of the corresponding Wolff potentials.  Theorem \ref{cor_main} brings
out the spectral side of this phenomenon, isolating the concentration
effect in the coefficient $\sum_{i=1}^{m_0} u_{n_0+i-1}^2(x_0)$
appearing in \eqref{eq:cor_main1}.

The rest of the paper is organized as follows. In Section
\ref{sec_ext} we introduce the well-known characterization of the
fractional Laplacian via the Caffarelli-Silvestre extension,
reformulating the definitions of $s$-fractional capacities within this
framework. In Section \ref{sec_pot}, we focus on the properties of
capacitary potentials.  In Section \ref{sec:appr-capac} we provide a
representation of the $s$-fractional $u$-capacity of $K$ in terms of
the integral of $u^2$ with respect to the $s$-equilibrium measure
$\mu_{\Omega,K}^s$, plus suitable remainder terms that are negligible
with respect to the capacity, see Proposition \ref{p:for-cap} and
Corollary \ref{cor:cor1}. This is the most technical part of the
paper, since the possible lack of regularity of $K$ and of the
corresponding capacitary potentials prevents direct integration by
parts; this difficulty is overcome by developing an approximation
procedure, which consists in approximating $K$ with smoother
$(N+1)$-dimensional compact sets and showing that the respective
capacitary potentials and capacities then approximate the original
ones.  In Section \ref{sec_eigen} we obtain an asymptotic expansion
for the eigenvalue variation and complete the proof of the main
theorems.

\section{An equivalent formulation by extension}\label{sec_ext}
Let $\R^{N+1}_+=\{z=(x,t)\in\R^{N+1}:\, x\in\R^N, t>0 \}$.
We define the
space $\mathcal{D}^{1,2}(\R^{N+1}_+;t^{1-2s})$ as the completion of
$C^\infty_c(\overline{\R^{N+1}_+})$ with respect to the norm
\begin{equation*}
\|U\|_{\mathcal{D}^{1,2}(\R^{N+1}_+;t^{1-2s})}= \sqrt{q(U,U)}
\end{equation*}
induced by the scalar product 
\begin{equation*}
q(U,V)=\int_{\R^{N+1}_+}t^{1-2s}\nabla U(x,t)\cdot\nabla V(x,t) \,dx\,dt.
\end{equation*}
We recall from \cite[Lemma 2.4]{FallFelli14} the following Hardy-type
inequality:
\begin{equation}\label{eq:hardy}
  \left(\frac{N-2s}{2}\right)^2 \int_{\R^{N+1}_+} t^{1-2s} \frac{U^2(z)}{|x|^2+t^2} \,dx\,dt \leq
  \int_{\R^{N+1}_+}t^{1-2s} | \nabla U(z)|^2 \,dx\,dt,
\end{equation}
for all $U\in \mathcal{D}^{1,2}(\R^{N+1}_+;t^{1-2s})$ and $s\in (0,\min\{1,N/2\})$. 

It is also well-known, see e.g. \cite{brandle2013concave}, that there
exists a well-defined continuous trace map
\begin{equation*}\label{eq:trace}
\mathop{\rm Tr}: \mathcal{D}^{1,2}(\R^{N+1}_+;t^{1-2s}) \to \mathcal{D}^{s,2}(\R^N)
\end{equation*}
which is onto. 
By the Caffarelli-Silvestre extension theorem
\cite{CaffarelliSilvestre07}, for every $u\in \mathcal D^{s,2}(\R^N)$, the minimum
\begin{equation*}
\min\left\{ q(W,W) : \, W\in \mathcal{D}^{1,2}(\R^{N+1}_+;t^{1-2s}),\, \mathop{\rm Tr}W=u\right\}
\end{equation*}
is achieved by a unique minimizer $U=\mathcal H(u) \in
\mathcal{D}^{1,2}(\R^{N+1}_+;t^{1-2s})$, which also 
satisfies 
\begin{equation}\label{eq:13}
q(\mathcal H(u),W)=\kappa_s(u,\mathop{\rm Tr}W)_{\mathcal
  D^{s,2}(\R^N)}\quad\text{for all }W \in
\mathcal{D}^{1,2}(\R^{N+1}_+;t^{1-2s}),
\end{equation}
where
\[
\kappa_s=\frac{\Gamma(1-s)}{2^{2s-1}\Gamma(s)},
\]
i.e. $U=\mathcal H(u)$
 weakly solves
\begin{equation*}
\begin{cases}
-\text{div}(t^{1-2s} \nabla U)=0, \quad&\text{in } \R^{N+1}_+, \\
\lim_{t\to0^+}\left(-t^{1-2s} \partial_t U \right)= \kappa_s
(-\Delta)^su, &\text{in } \R^N\times\{0\}.
\end{cases}
\end{equation*}
We observe that \eqref{eq:13} implies that
\begin{equation}\label{eq:4}
  \|U\|_{\mathcal{D}^{1,2}(\R^{N+1}_+;t^{1-2s})}^2=\kappa_s\|u\|_{\mathcal D^{s,2}(\R^N)}^2.
\end{equation}
Therefore, if $\lambda_j^s(\Omega)$ is an eigenvalue of
\eqref{eq:eigenvalue_main} for a certain $j\in\N\setminus\{0\}$
and $u_j \in \mathcal D^{s,2}(\Omega)$ is an associated eigenfunction,
the extension $U_j=\mathcal H(u_j)$ satisfies $\mathop{\rm Tr} U_j=u_j$ and
\begin{equation}\label{eq:eige-exte}
\begin{cases}
-\text{div}(t^{1-2s} \nabla U_j)=0, \quad&\text{in } \R^{N+1}_+, \\
\lim_{t\to0^+}\left(-t^{1-2s} \partial_t U_j
\right)=\lambda_j^s(\Omega) \kappa_s \mathop{\rm Tr} U_j,
 &\text{in } \Omega\times\{0\}, \\
U_j=0, &\text{in } (\R^N\setminus \Omega)\times\{0\},
\end{cases}
\end{equation}
in a weak sense, that is
\begin{equation*}
\begin{cases}
U_j\in
  \mathcal{E}^{s}(\Omega),\\
{\displaystyle{ q(U_j,\phi)
    =\lambda_j^s(\Omega) \kappa_s\int_\Omega \mathop{\rm Tr}U_j
    \mathop{\rm Tr}\phi\,dx}}\quad\text{for every }\phi\in
\mathcal{E}^{s}(\Omega),
\end{cases}
\end{equation*}
where $\mathcal{E}^{s}(\Omega)$ is defined as 
the
closure of $C^\infty_c(\R^{N+1}_+\cup\Omega)$ in
$\mathcal{D}^{1,2}(\R^{N+1}_+;t^{1-2s})$, equivalently characterized
as 
\begin{equation*}
\mathcal{E}^{s}(\Omega)
=\{U\in \mathcal{D}^{1,2}(\R^{N+1}_+;t^{1-2s}): \, 
\mathop{\rm Tr}U\in \mathcal D^{s,2}(\Omega)\}.
\end{equation*}
In particular, in view of \eqref{eq:4}, the
$s$-fractional capacities introduced in Definition
\ref{def:s-capacity}  can be equivalently characterized as 
\begin{equation}\label{eq:def_capa}
\mathop{\rm Cap}\nolimits^s_\Omega(K)=\frac1{\kappa_s}\inf \left\{
q(W,W): \, W\in \mathcal{E}^{s}(\Omega) \text{ and } 
W-\eta_K \in \mathcal{E}^{s}(\Omega\setminus K)\right\},
\end{equation}
$\eta_K \in C^\infty_c(\R^{N+1}_+\cup\Omega)$ being any fixed
function such that $\eta_K=1$ in a neighborhood of $K$, and, for every
$u\in \mathcal D^{s,2}(\Omega)$, 
\begin{equation}\label{eq:capa-u}
\mathop{\rm Cap}\nolimits^s_\Omega(K,u)=\frac1{\kappa_s}\inf \left\{
q(W,W):  W\in \mathcal{E}^{s}(\Omega),\,
W-U \in \mathcal{E}^{s}(\Omega\setminus K)
\right\},
\end{equation}
where $U\in  \mathcal{E}^{s}(\Omega)$ is
such that $\mathop{\rm Tr}U=u$.

\section{Capacitary potentials}\label{sec_pot}
By classical minimization methods, the infimum in \eqref{eq:def_capa}
is realized by a unique minimizer, i.e., there exists a unique
$V_{\Omega,K} \in \mathcal{E}^{s}(\Omega)$ such that
$V_{\Omega,K}-\eta_K \in \mathcal{E}^{s}(\Omega\setminus K)$ and
\begin{equation*}
\mathop{\rm Cap}\nolimits^s_\Omega(K)=\frac1{\kappa_s} q(V_{\Omega,K} ,V_{\Omega,K}).
\end{equation*}
Such a capacitary potential  $V_{\Omega,K}$ is the unique weak solution to 
\begin{equation*}
\begin{cases}
-\text{div}(t^{1-2s}\nabla V_{\Omega,K})=0, \quad&\text{in } \R^{N+1}_+, \\
\lim_{t\to0^+}\left(-t^{1-2s} \partial_t V_{\Omega,K} \right)=0, &\text{in }
(\Omega\setminus K)\times\{0\}, \\
V_{\Omega,K}=0, &\text{in } (\R^N\setminus\Omega)\times\{0\}, \\
V_{\Omega,K}=1, &\text{in } K\times\{0\},
\end{cases}
\end{equation*}
i.e.  $V_{\Omega,K}\in \mathcal E^s(\Omega)$  is the unique function such that 
\begin{equation}\label{eq:eVK}
V_{\Omega,K}-\eta_K \in \mathcal E^s(\Omega\setminus K)\quad\text{and}
\quad 
q(V_{\Omega,K},\phi )=0 \quad\text{for every } \phi\in
\mathcal E^s(\Omega\setminus K).
\end{equation}
We observe that
  \begin{equation*}
    \mathop{\rm Tr}V_{\Omega,K}= v_{\Omega,K},
  \end{equation*}
where $v_{\Omega,K}$ is the unique minimizer in \eqref{eq:def-capOK},
see Section \ref{subsec_equil_measu}.

Similarly,  for every
$u\in \mathcal D^{s,2}(\Omega)$,  the
infimum in \eqref{eq:capa-u}  is achieved by a unique capacitary potential
$V_{\Omega,K,u}\in \mathcal{E}^s(\Omega)$, so
that
\begin{equation}\label{eq_capa_VKu}
\mathop{\rm Cap}\nolimits^s_\Omega(K,u)= \frac1{\kappa_s}q(
V_{\Omega,K,u}, V_{\Omega,K,u})
\end{equation}
and $V_{\Omega,K,u}$ is the unique weak solution of
\begin{equation*}
\begin{cases}
-\text{div}(t^{1-2s} \nabla V_{\Omega,K,u})=0, \quad&\text{in } \R^{N+1}_+, \\
\lim_{t\to0^+}\left(-t^{1-2s} \partial_t V_{\Omega,K,u} \right)=0,
&\text{in }(\Omega\setminus K)\times\{0\},\\
 V_{\Omega,K,u}=0, &\text{in } (\R^N\setminus\Omega)\times\{0\}, \\
V_{\Omega,K,u}=u, &\text{in } K\times\{0\},
\end{cases}
\end{equation*}
in the sense that $V_{\Omega,K,u}\in \mathcal{E}^s(\Omega)$,
$u-\mathop{\rm Tr}V_{\Omega,K,u}\in \mathcal D^{s,2}(\Omega\setminus
K)$, and
\begin{equation}\label{eq:VKu}
q(V_{\Omega,K,u},\phi )=0 \quad\text{for all } \phi\in 
\mathcal{E}^s(\Omega\setminus K).
\end{equation}
We also define
\begin{equation}\label{eq:pot-tr}
  v_{\Omega,K}^u:=\mathop{\rm Tr}V_{\Omega,K,u}
\end{equation}
and observe that $v_{\Omega,K}^u\in \mathcal D^{s,2}(\Omega)$ is the
unique weak solution to
\begin{equation*}
\begin{cases}
(-\Delta)^s v_{\Omega,K}^u=0, \quad& \text{in } \Omega\setminus K, \\
v_{\Omega,K}^u=0, &\text{in }\R^N\setminus\Omega,\\
v_{\Omega,K}^u=u, &\text{in }K,
\end{cases}
\end{equation*}
i.e. $v_{\Omega,K}^u-u\in \mathcal D^{s,2}(\Omega\setminus K)$ and 
\begin{equation}\label{eq_vKu}
  (v_{\Omega,K}^u,\phi)_{\mathcal D^{s,2}(\R^N)}=0\quad\text{for
    every }\phi \in \mathcal D^{s,2}(\Omega\setminus K).
\end{equation}
Furthermore, by \eqref{eq:4} and \eqref{eq_capa_VKu}
\begin{equation}\label{eq_capa_vKu}
  \mathop{\rm Cap}\nolimits^s_\Omega(K,u)=
  \frac1{\kappa_s} q(V_{\Omega,K,u} ,V_{\Omega,K,u})=\|v_{\Omega,K}^u\|^2_{\mathcal{D}^{s,2}(\Omega)}.
\end{equation}
The following characterization of the capacities
$\mathop{\rm Cap}\nolimits^s_\Omega(K,u)$ and
$\mathop{\rm Cap}\nolimits^s_\Omega(K)$ follows from density arguments
and will be useful in the sequel.
\begin{lemma}\label{l:cara-cap}
  If $u\in C^\infty(\Omega)$ and $K\subset\Omega$ is
  a compact set, then
  \begin{equation*}
    \mathop{\rm Cap}\nolimits^s_\Omega(K,u)=
\frac1{\kappa_s}\inf_{W\in \mathcal D^u_{\Omega,K}}
q(W,W)
\end{equation*}
where
\begin{equation*}
  \mathcal D^u_{\Omega,K}=\left\{
  \begin{array}{ll}
    W\in \mathcal{E}^{s}(\Omega)\cap
C^0(\overline{\R^{N+1}_+}):
      \hskip-5pt&W(x,t)-u(x)=0\text{ for every $(x,t)\in \mathcal U\cap
  \overline{\R^{N+1}_+}$} \\[3pt]
    &\text{for some $\mathcal U$ open set in $\R^{N+1}$
      with $K\times\{0\}\subset \mathcal U$}
  \end{array}
  \right\}.
      \end{equation*}
      In particular
      \begin{equation*}
    \mathop{\rm Cap}\nolimits^s_\Omega(K)=
\frac1{\kappa_s}\inf_{W\in \mathcal D^1_{\Omega,K}}
q(W,W).
\end{equation*}
\end{lemma}
\begin{proof}
  Let $\zeta_K\in C^\infty_c(\Omega)$ be such that $\zeta_K=1$ in a
  neighborhood of $K$. Since $V_{\Omega,K,u}\in \mathcal{E}^s(\Omega)$
  and
  $\mathop{\rm Tr}V_{\Omega,K,u}-\zeta_K u\in \mathcal
  D^{s,2}(\Omega\setminus K)$, there exists a sequence of smooth
  functions $\{v_n\}\subset C^\infty_{\rm c}(\Omega\setminus K)$ such
  that $v_n\to \mathop{\rm Tr}V_{\Omega,K,u}-\zeta_K u$ in
  $\mathcal D^{s,2}(\Omega\setminus K)$. Then, denoting as
  $V_n:=\mathcal H(v_n)$, we have
  \begin{equation*}
  V_n\to V_{\Omega,K,u}-\mathcal H(\zeta_Ku) \quad\text{in }
  \mathcal D^{1,2}(\R^{N+1}_+;t^{1-2s}).
\end{equation*}
We observe that $V_n+\mathcal H(\zeta_Ku)\in
C^0(\overline{\R^{N+1}_+})$ with $\mathop{\rm Tr}(
V_n+\mathcal H(\zeta_Ku))=u$ in an $N$-dimensional neighborhood
of $K$.  Furthermore, $q(V_n+\mathcal H(\zeta_Ku), V_n+\mathcal
H(\zeta_Ku))\to \kappa_s\mathop{\rm Cap}\nolimits^s_\Omega(K,u)$.
Hence
 \begin{equation*}
    \mathop{\rm Cap}\nolimits^s_\Omega(K,u)=
\frac1{\kappa_s}\inf_{W\in \widetilde{\mathcal D}^u_{\Omega,K}}
q(W,W)
\end{equation*}
where
\begin{equation*}
  \widetilde{\mathcal D}^u_{\Omega,K}=\left\{
  \begin{array}{ll}
    W\in \mathcal{E}^{s}(\Omega)\cap
C^0(\overline{\R^{N+1}_+}):
      \hskip-5pt&W(x,0)- u(x)=0\text{ for every $x\in U$} \\[3pt]
    &\text{for some $U$ open set in $\R^{N}$
      with $K\subset U$}
  \end{array}
  \right\}.
      \end{equation*}
Let us fix some 
  $\eta\in C^\infty_{\rm c}([0,+\infty))$ such that $\eta(t)=1$ for
  all $t\in [0,1]$ and $G\in C^\infty(\R)$ such that $G(t)=0$ for
  $t\in [-1,1]$ and $G(t)=t$ for $|t|\geq2$.
 The conclusion follows by observing that, for every $W\in
\widetilde{\mathcal D}^u_{\Omega,K}$,
\begin{equation*}
  W_n(x,t):=\frac1n G\big(n\big(W(x,t)-\eta(t)\zeta_K(x) u(x)\big)\big)+\eta(t)\zeta_K(x) u(x)\to W\quad\text{in
  }\mathcal D^{1,2}(\R^{N+1}_+;t^{1-2s})
\end{equation*}
and $W_n\in \mathcal D^u_{\Omega,K}$.
\end{proof}

The following lemma provides a slight generalization of \cite[Lemma
4.1]{AFN}.
\begin{lemma}\label{lem:trace_o_capacity}
  Let $\{K(\e)\}_{\e>0}$ be a family of compact sets contained in the
  bounded open set $\Omega\subset\R^N$ such that
  $\lim_{\e\to 0^+}\mathop{\rm Cap}\nolimits^s_\Omega(K(\e))=0$. Then
\begin{equation}\label{eq:o1}
\int_\Omega |\mathop{\rm Tr}V_{\Omega,K(\e)}|^2\,dx 
=o(\mathop{\rm Cap}\nolimits^s_\Omega(K(\e)))
\quad\text{as } \varepsilon\to 0,
\end{equation}
and, for every bounded open set
  $\mathcal U\subset\R^{N+1}_+$,
\begin{equation}\label{eq:o2}
  \int_{\mathcal U} t^{1-2s}V_{\Omega,K(\e)}^2\,dx\,dt 
=o(\mathop{\rm Cap}\nolimits^s_\Omega(K(\e)))
\quad\text{as } \varepsilon\to 0.
\end{equation}
\end{lemma}
\begin{proof}
  To prove \eqref{eq:o1}, we argue by contradiction and assume that
  there exist a sequence $\e_n\to0^+$ and a
  constant $C>0$ such that
\begin{equation*}
\int_\Omega |\mathop{\rm Tr}V_{\Omega,K(\e_n)}|^2\,dx 
\geq C
\int_{\R^{N+1}_+}t^{1-2s}|\nabla V_{\Omega,K(\e_n)} |^2 \,dx\,dt
\end{equation*}
for every $n$. Letting
\begin{equation*}
W_n=\frac{V_{\Omega,K(\e_n)}}{\| \mathop{\rm Tr} V_{\Omega,K(\e_n)}\|_{L^2(\Omega)}},
\end{equation*}
we have
\begin{equation*}
\|\mathop{\rm Tr} W_n\|_{L^2(\Omega)}=1 \qquad\text{and}\qquad
\|W_n\|_{\mathcal{D}^{1,2}(\R^{N+1}_+;t^{1-2s})}^2\leq C^{-1}
\end{equation*}
for every $n$. Then, 
there exist a subsequence $\{n_k\}_{k\geq1}$ and $W\in \mathcal
E^s(\Omega)$ such that
\begin{equation*}
W_{n_k} \rightharpoonup W \quad \text{in }\mathcal{D}^{1,2}(\R^{N+1}_+;t^{1-2s}) \text{ as }k\to+\infty 
\end{equation*}
and, due to compactness
of the trace operator $\mathop{\rm
  Tr}:\mathcal E^s(\Omega) \to
L^2(\Omega)$,
\begin{equation}\label{eq:trw0}
\|\mathop{\rm Tr} W\|_{L^2(\Omega)}=1.
\end{equation}
Let $\phi\in C^\infty_c(\R^{N+1}_+\cup\Omega)$. Then,
$\phi(1-V_{\Omega,K(\e_n)})\in \mathcal E^s(\Omega\setminus K(\e_n))$, so
that, by \eqref{eq:eVK},
\begin{align}
\label{eq:ptl}  0&=\int_{\R^{N+1}_+}t^{1-2s}\nabla W_n\cdot \nabla \phi\,dx\,dt\\
 \notag &\quad -
  \int_{\R^{N+1}_+}t^{1-2s}\nabla W_n\cdot \nabla \phi\, V_{\Omega,K(\e_n)}\,dx\,dt-
  \int_{\R^{N+1}_+}t^{1-2s}\nabla W_n\cdot \nabla
  V_{\Omega,K(\e_n)}\phi \,dx\,dt.
\end{align}
By \eqref{eq:hardy},
\begin{align*}
  &\left| \int_{\R^{N+1}_+}t^{1-2s}\nabla W_n\cdot \nabla \phi\,
    V_{\Omega,K(\e_n)}\,dx\,dt\right|\\
  &  \leq \Big(\sup_{\R^{N+1}_+}|\nabla
    \phi|\Big) \Big(\sup_{z\in\mathop{\rm
    supp\phi}}|z|\Big)\sqrt{\int_{\R^{N+1}_+}
    t^{1-2s}|\nabla W_n|^2\,dx\,dt}
    \sqrt{\int_{\R^{N+1}_+}t^{1-2s}
    \frac{V_{\Omega,K(\e_n)}^2}{|x|^2+t^2}\,dx\,dt}\\
  &\leq \mathop{\rm const}\sqrt{\mathop{\rm Cap}\nolimits^s_\Omega(K(\e_n))}
\end{align*}
for some $\mathop{\rm const}>0$ independent of $n$. In a similar way,
\begin{align*}
&\left| \int_{\R^{N+1}_+}t^{1-2s}\nabla W_n\cdot \nabla V_{\Omega,K(\e_n)}\phi \,dx\,dt\right|  \\
&\qquad\leq  \sqrt{\kappa_s}\|\phi\|_{L^\infty(\R^{N+1}_+)}
\sqrt{\mathop{\rm Cap}\nolimits^s_\Omega(K(\e_n))}\sqrt{\int_{\R^{N+1}_+}t^{1-2s}|\nabla W_n|^2\,dx\,dt}\\
&\qquad\leq   \mathop{\rm const}\sqrt{\mathop{\rm Cap}\nolimits^s_\Omega(K(\e_n))}.
\end{align*}
Hence, passing to the limit in \eqref{eq:ptl} as $n=n_k\to\infty$
yields
\begin{equation*}
  \int_{\R^{N+1}_+}t^{1-2s}\nabla W\cdot \nabla \phi\,dx\,dt=0
\end{equation*}
for every $\phi\in C^\infty_c(\R^{N+1}_+\cup\Omega)$ and hence, by
density, for every $\phi\in\mathcal E^s(\Omega)$. This implies that
$W\equiv0$, thus contradicting \eqref{eq:trw0}.

To prove \eqref{eq:o2}, we can argue in a similar way, exploiting the
compactness of the restriction map from
$\mathcal{D}^{1,2}(\R^{N+1}_+;t^{1-2s})$ into
$L^2(\mathcal U;t^{1-2s})$.
\end{proof}

\section{Approximating capacities}\label{sec:appr-capac}

The possible lack of regularity of the hole (which is merely assumed
to be a compact set) and of its capacitary potentials creates several
technical difficulties, such as preventing direct integration by parts
in integrals involving the potentials. Hence, an approximation
strategy is developed in this section: this consists of approximating
the hole by a family of smoother $(N+1)$-dimensional compact sets and
showing the convergence of the corresponding capacitary potentials and
capacities to the original ones.

 Let $\Omega\subset\R^N$ be a bounded open set and
 $s\in (0,\min\{1,N/2\})$.

 Let $\mathcal K$ be a compact set in
$\overline{\R^{N+1}_+}$ such that
\begin{align}
\label{eq:calK1} &\text{$\widetilde{\mathcal
  K}=\{(x,t)\in \R^{N+1}:(x,|t|)\in\mathcal K\}$ is the closure of a
$C^\infty$-smooth open set in $\R^{N+1}$};\\
  \label{eq:calK3}&K=\{(x,t)\in\mathcal K:t=0\}\subset
                    \Omega\times\{0\}\quad\text{and}
                    \quad \mathcal K\subseteq K\times[0,1); \\
\label{eq:calK2} &\mathcal K\cap
\{(x,t)\in\R^{N+1}:0\leq t\leq \delta_0\}=K\times[0,\delta_0] \text{ for some
$\delta_0>0$}.
\end{align}
 For any function $u\in C^\infty(\Omega)$, we introduce
 the following notions of  capacity and $u$-capacity of $\mathcal K$:
\begin{equation}\label{eq:app1}
\mathcal C^s_\Omega(\mathcal K)=\frac1{\kappa_s}\inf \left\{
q(W,W): \, W\in \mathcal{E}^{s}(\Omega) \text{ and } 
W(x,t)-1=0\text{ a.e. in }\mathcal K\right\},
\end{equation}
and
\begin{equation}\label{eq:app2}
  \mathcal C^s_\Omega(\mathcal K,u)
=\frac1{\kappa_s}\inf \left\{
  q(W,W):  W\in \mathcal{E}^{s}(\Omega)
\text{ and } 
W(x,t)-u(x)=0\text{ a.e. in }\mathcal K\right\}.
\end{equation}
We denote by $\mathcal V_{\mathcal K}$ and $\mathcal V_{\mathcal K,u}$
the unique functions attaining
$\mathcal C^s_\Omega(\mathcal K)$ and $\mathcal C^s_\Omega(\mathcal
K,u)$, respectively. In particular,  $\mathcal V_{\mathcal K}$ and $\mathcal V_{\mathcal K,u}$
satisfy
\begin{equation}\label{eq:calVKu}
  q(\mathcal V_{\mathcal K},\phi )=q(\mathcal V_{\mathcal K,u},\phi )=0 \quad\text{for all } \phi\in 
\mathcal{E}^s(\Omega)\text{ such that $\phi=0$ a.e. in $\mathcal K$}.
\end{equation}
By classical elliptic regularity, we have that 
\begin{equation*}
  \mathcal V_{\mathcal K},\mathcal V_{\mathcal K,u}\in
  C^\infty(\R^{N+1}_+\setminus \mathring{\mathcal K}),
\end{equation*}
where $\mathring{\mathcal K}$ denotes the interior of $\mathcal K$.
Further properties of the functions $\mathcal V_{\mathcal K},\mathcal
V_{\mathcal K,u}$ are established in the following lemma.
\begin{lemma}\label{l:prop_calV}
  Let $\Omega\subset\R^N$ be a bounded open set and
  $s\in (0,\min\{1,N/2\})$.  Let
  $\mathcal K$ be a compact set in
  $\overline{\R^{N+1}_+}$ satisfying \eqref{eq:calK1}, \eqref{eq:calK3}, and
  \eqref{eq:calK2}.  For some $u\in C^\infty(\Omega)$, let
  $\mathcal V_{\mathcal K}$ and $\mathcal V_{\mathcal K,u}$ be the
  unique functions attaining $\mathcal C^s_\Omega(\mathcal K)$ and
  $\mathcal C^s_\Omega(\mathcal K,u)$, respectively. Then the
  following properties hold:
  \begin{enumerate}[\rm (i)]
\item $\mathcal V_{\mathcal K}\in
  C^{1}((\Omega\setminus \mathring{K})\times[0,\delta_0))$;
\item
  $t^{1-2s}\partial_t \mathcal V_{\mathcal K}\in
  C^{0}((\Omega\setminus K)\times[0,\delta_0))$ and
  $\lim_{t\to0^+}t^{1-2s}\partial_t \mathcal V_{\mathcal K}(x,t)=0$
  for all $x\in \Omega\setminus K$;
\item $0\leq \mathcal V_{\mathcal K}(x,t)\leq 1$ for all $(x,t)\in
  \R^{N+1}_+$;
\item
  $|\mathcal V_{\mathcal K,u}(x,t)|\leq \left(\max_{\mathcal
      K}|u|\right) \mathcal V_{\mathcal K}(x,t)$ for all
  $(x,t)\in \R^{N+1}_+$;
\item $\partial_\nu \mathcal V_{\mathcal K}\geq 0$ for every
  $(x,t)\in \partial \mathcal K\cap\R^{N+1}_+$, $\nu$ being the unit vector
  on $\partial\mathcal K$ pointing toward the interior of
  $\mathcal K$;
\item  $\mathcal C^s_\Omega(\mathcal K)=\frac{1}{\kappa_s}\int_{\partial \mathcal
    K\cap \R^{N+1}_+} t^{1-2s}\partial_\nu \mathcal V_{\mathcal
    K}\,dS=\frac{1}{\kappa_s}
  \int_{\partial \mathcal
    K\cap \R^{N+1}_+} t^{1-2s}|\partial_\nu \mathcal V_{\mathcal K}|\,dS$.
\end{enumerate}

\end{lemma}
\begin{proof}
  Statements (i) and (ii) directly follow from \cite[Lemma A.1]{dLFV2022},
  taking into account assumption \eqref{eq:calK2}.

To prove (iii), we observe that $(\mathcal V_{\mathcal K}-1)^+\in 
\mathcal{E}^s(\Omega)$ and $(\mathcal V_{\mathcal K}-1)^+=0$ a.e. in
$\mathcal K$. Hence we can choose $\phi= (\mathcal V_{\mathcal
  K}-1)^+$ in \eqref{eq:calVKu}, thus obtaining 
\begin{equation*}
0=\int_{\R^{N+1}_+}t^{1-2s}\nabla \mathcal V_{\mathcal K} \cdot\nabla
(\mathcal V_{\mathcal K}-1)^+ \,dx\,dt=
\int_{\R^{N+1}_+}t^{1-2s}|\nabla (\mathcal V_{\mathcal K}-1)^+|^2 \,dx\,dt  ,
\end{equation*}
and, consequently, $(\mathcal V_{\mathcal K}-1)^+\equiv 0$ in
$\R^{N+1}_+$, i.e. $\mathcal V_{\mathcal K}\leq 1$ in $\R^{N+1}_+$.
Similarly, 
$\mathcal V_{\mathcal K}^-\in 
\mathcal{E}^s(\Omega)$ and $\mathcal V_{\mathcal K}^-=0$ in
$\mathcal K$, so that, choosing  $\phi= \mathcal V_{\mathcal
  K}^-$ in \eqref{eq:calVKu}, we obtain $\mathcal V_{\mathcal K}^-\equiv 0$ in
$\R^{N+1}_+$, i.e. $\mathcal V_{\mathcal K}\geq 0$ in $\R^{N+1}_+$.
Statement (iv) follows in an analogous way by choosing in \eqref{eq:calVKu}
first $\phi= \big(\mathcal V_{\mathcal
  K,u}-(\max_{\mathcal K}|u|)\mathcal V_K\big)^+$, and then 
$\phi= \big(\mathcal V_{\mathcal
  K,u}+(\max_{\mathcal K}|u|)\mathcal V_K\big)^-$.

Statement (v) directly follows from (iii). Finally, (vi) can be
obtained by an integration by parts, taking into account (v) and (ii).
\end{proof}

In the following lemma, we observe that
$\mathop{\rm Cap}\nolimits^s_\Omega(K,u)$ and
$\mathop{\rm Cap}\nolimits^s_\Omega(K)$ can be approximated by 
capacities of type \eqref{eq:app1}-\eqref{eq:app2}.

\begin{lemma}\label{l:appr-cap_n}
  Let $\Omega\subset\R^N$ be a bounded open set and
  $s\in (0,\min\{1,N/2\})$.  Let $\{\mathcal K_n\}_n$ be a sequence of
  compact sets in $\overline{\R^{N+1}_+}$ such that $\mathcal K_n$
  satisfies \eqref{eq:calK1}, \eqref{eq:calK3}, and \eqref{eq:calK2}
  for every $n$ and 
 \begin{equation*}
  \mathcal K_{n+1}\subseteq \mathcal K_{n}\text{ for all $n$},\quad 
  \bigcap_{n}\mathcal K_n=K\times\{0\},
\end{equation*}
for some  compact set $K\subset\Omega$.   Then
    \begin{equation}\label{eq:approx1}
     \lim_{n\to\infty}  \mathcal C^s_\Omega(\mathcal K_n)
     =\mathop{\rm Cap}\nolimits^s_\Omega(K),
   \end{equation}
and,    for every $u\in C^\infty(\Omega)$,
   \begin{equation}\label{eq:approx2}
     \lim_{n\to\infty}  \mathcal C^s_\Omega(\mathcal K_n,u)
     =\mathop{\rm Cap}\nolimits^s_\Omega(K,u).
   \end{equation}
   Furthermore,
   \begin{equation}\label{eq:con-pot}
     \mathcal V_{\mathcal K_n,u}\to V_{\Omega,K,u}\quad\text{and}\quad
          \mathcal V_{\mathcal K_n}\to V_{\Omega,K}
        \end{equation}
        in $\mathcal{D}^{1,2}(\R^{N+1}_+;t^{1-2s})$.
 \end{lemma}
\begin{proof}
  We first focus on the proof of \eqref{eq:approx2}, which includes
  \eqref{eq:approx1} as a particular case.
Since the sequence $\{\mathcal C^s_\Omega(\mathcal K_n,u)\}_n$ is
decreasing, there exists a limit
\begin{equation*}
L:=  \lim_{n\to\infty}\mathcal C^s_\Omega(\mathcal K_n,u)=\inf_n \mathcal
  C^s_\Omega(\mathcal K_n,u)\in [0,+\infty).
\end{equation*}
  Since $K\times\{0\}\subset\mathcal
  K_n$ for all $n$, the set of values which is minimized in
  \eqref{eq:app2} with $\mathcal K=\mathcal K_n$ is contained in the
  set  which is minimized in
  \eqref{eq:capa-u} for $K$, so that
  \begin{equation*}
    \mathop{\rm Cap}\nolimits^s_\Omega(K,u)\leq   \mathcal C^s_\Omega(\mathcal K_n,u)
 \quad\text{for every }n.
\end{equation*}
Hence
  \begin{equation*}
    \mathop{\rm Cap}\nolimits^s_\Omega(K,u)\leq   L.
\end{equation*}
On the other hand, in view of Lemma \ref{l:cara-cap}, for every
$\e>0$, there exists $W\in \mathcal D^u_{\Omega,K}$ such that
\begin{equation*}
  \frac1{\kappa_s}q(W,W)<\mathop{\rm Cap}\nolimits^s_\Omega(K,u)+\e.
\end{equation*}
Since $W\in \mathcal D^u_{\Omega,K}$, there exists an open set
$\mathcal U\subset \R^{N+1}$ such that
$K\times\{0\}\subset \mathcal U$ and $W(x,t)-u(x)=0$ for every
$(x,t)\in \mathcal U\cap \overline{\R^{N+1}_+}$. We observe that there
exists $n_0$ such that
$\mathcal K_{n_0}\subset \mathcal U\cap
\overline{\R^{N+1}_+}$. Therefore, $W(x,t)-u(x)=0$ for every
$(x,t)\in \mathcal K_{n_0}$, hence
  \begin{equation*}
  L=\inf_n \mathcal
  C^s_\Omega(\mathcal K_n,u)\leq
  \mathcal C^s_\Omega(\mathcal K_{n_0},u)\leq
    \frac1{\kappa_s}q(W,W)
    <\mathop{\rm Cap}\nolimits^s_\Omega(K,u)+\e.
  \end{equation*}
  Hence $L< \mathop{\rm Cap}\nolimits^s_\Omega(K,u)+\e$ for every
  $\e>0$, so that necessarily
  $L\leq \mathop{\rm Cap}\nolimits^s_\Omega(K,u)$ and the proof of
  \eqref{eq:approx2} is complete.
  
  To prove \eqref{eq:con-pot}, we observe that
  \begin{equation*}
\mathcal V_{\mathcal K_n,u}- V_{\Omega,K,u}\in \mathcal E^s(\Omega\setminus K),   
\end{equation*}
therefore we can test \eqref{eq:VKu} with $\mathcal V_{\mathcal K_n,u}-
V_{\Omega,K,u}$ thus obtaining
\begin{equation*}
  \int_{\R^{N+1}_+}t^{1-2s}|\nabla V_{\Omega,K,u}|^2 \,dx\,dt=
   \int_{\R^{N+1}_+}t^{1-2s}\nabla V_{\Omega,K,u}\cdot\nabla  \mathcal V_{\mathcal K_n,u}\,dx\,dt.
 \end{equation*}
 Hence,
 \begin{multline*}
   \int_{\R^{N+1}_+}t^{1-2s}|\nabla( V_{\Omega,K,u}-\mathcal
   V_{\mathcal K_n,u})|^2 \,dx\,dt=
   \int_{\R^{N+1}_+}t^{1-2s}|\nabla V_{\Omega,K,u}|^2 \,dx\,dt\\+
   \int_{\R^{N+1}_+}t^{1-2s}|\nabla \mathcal
   V_{\mathcal K_n,u}|^2 \,dx\,dt-2 \int_{\R^{N+1}_+}t^{1-2s}|\nabla
   V_{\Omega,K,u}|^2 \,dx\,dt\to 0
 \end{multline*}
 as $n\to\infty$ in view of \eqref{eq:approx2}. The proof is thereby complete. 
\end{proof}

\begin{remark}
  If $K\subset \Omega$ is a compact set in the bounded open set
  $\Omega\subset\R^N$, we can consider a sequence $\{\mathcal K_n\}_n$
  of compact sets in $\overline{\R^{N+1}_+}$ such that, for every $n$,
  $\mathcal K_n$ satisfies \eqref{eq:calK1}, \eqref{eq:calK3}, and
  \eqref{eq:calK2}, $\mathcal K_{n+1}\subseteq \mathcal K_{n}$, and
  $\bigcap_{n}\mathcal K_n=K \times\{0\}$. Then, from Lemma
  \ref{l:prop_calV}-(iii), (iv) and Lemma \ref{l:appr-cap_n} it follows
  that
  \begin{equation*}
  0\leq V_{\Omega, K}(x,t)\leq 1\quad\text{for all }(x,t)\in
  \R^{N+1}_+,
\end{equation*}
and, for every $u\in C^\infty(\Omega)$,
  \begin{equation}\label{eq:rs}
| V_{\Omega,K,u}(x,t)|\leq \left(\max_K|u|\right)
V_{\Omega,K}(x,t)\quad\text{for all }
(x,t)\in \R^{N+1}_+.
\end{equation}
\end{remark}
The potentials $\mathcal V_{\mathcal K,u}$ and
$\mathcal V_{\mathcal K}$ enjoy enough regularity to allow integration
by parts once \eqref{eq:calVKu} is tested with
$\mathcal V_{\mathcal K,u}-\mathcal V_{\mathcal K} u$, thus allowing
the capacity $\mathcal C^s_\Omega(\mathcal K,u)$ to be written as an
integral involving the normal derivative of $\mathcal V_{\mathcal K}$
on the boundary of $\mathcal K$ plus some remainder terms.

\begin{lemma}\label{l:stimecap}
  Let $\Omega\subset\R^N$ be a bounded open set and
  $s\in (0,\min\{1,N/2\})$.  Let
  $\mathcal K\subset \overline{\R^{N+1}_+}$ be a compact set in
  $\overline{\R^{N+1}_+}$ satisfying \eqref{eq:calK1},
  \eqref{eq:calK3}, and \eqref{eq:calK2}.  Let
  $\eta\in C^\infty_{\rm c}([0,+\infty))$ be such that $\eta(t)=1$ for
  all $t\in [0,1]$ and $\zeta_K\in C^\infty_c(\Omega)$ be such that
  $\zeta_K=1$ in a neighborhood of $K$, where
  $K=\{(x,t)\in\mathcal K:t=0\}\subset\Omega\times\{0\}$. If
  $u\in C^\infty(\Omega)$, then
  \begin{align*}
\mathcal C^s_\Omega(\mathcal K,u)=&
\frac{1}{\kappa_s}\int_{\partial\mathcal K\cap\R^{N+1}_+}t^{1-2s}u^2\partial_\nu
    \mathcal V_{\mathcal K}\,dS- \frac{2}{\kappa_s} \int_{\R^{N+1}_+\setminus\mathcal
                                     K}t^{1-2s}\mathcal V_{\mathcal K,u}\nabla
    \mathcal V_{\mathcal K}\cdot\nabla (u\eta\zeta_K)\,dx\,dt\\
&\quad    -\frac{1}{\kappa_s}\int_{\R^{N+1}_+\setminus\mathcal K}\mathop{\rm
                                                                  div}(t^{1-2s}\nabla
                                                                  (u\eta\zeta_K))\mathcal
                                                                  V_{\mathcal
                                                                  K,u}
    \mathcal V_{\mathcal K}\,dx\,dt-\frac{1}{\kappa_s}\int_{\mathcal K}t^{1-2s}u\Delta u\,dx\,dt,
  \end{align*}
  where $\nu$ is the unit vector on $\partial\mathcal K$ pointing
  toward the interior of $\mathcal K$.
\end{lemma}
\begin{proof}
We observe that, by \eqref{eq:calK3}, $\eta(t)\zeta_K(x)= 1$ for every
$(x,t)\in \mathcal K$.  Then we can test \eqref{eq:calVKu} with
  $\phi=\mathcal V_{\mathcal K,u}-\mathcal V_{\mathcal K} u\eta\zeta_K$, thus obtaining
\begin{align*}
  \int_{\R^{N+1}_+\setminus\mathcal K}&t^{1-2s}|\nabla \mathcal V_{\mathcal K,u}|^2\,dx\,dt=
  \int_{\R^{N+1}_+\setminus\mathcal K}t^{1-2s}\nabla \mathcal V_{\mathcal K,u}
  \cdot \nabla (\mathcal V_{\mathcal K} u\eta\zeta_K)\,dx\,dt\\
                                      &=\lim_{\delta\to0^+}
                                        \int_{\{(x,t)\in\R^{N+1}_+:t>\delta\}\setminus\mathcal
                                        K}
                                        t^{1-2s}\nabla \mathcal
                                        V_{\mathcal K,u} \cdot \nabla (\mathcal V_{\mathcal K} u\eta\zeta_K)\,dx\,dt\\
 \\
                                      &= \lim_{\delta\to0^+}\bigg(
                                        -2\int_{\{(x,t)\in\R^{N+1}_+:t>\delta\}\setminus\mathcal
                                        K}
                                        t^{1-2s}\mathcal
                                        V_{\mathcal K,u}\nabla \mathcal V_{\mathcal K}
                                        \cdot\nabla
                                        (u\eta\zeta_K)\,dx\,dt\\
&\qquad\qquad -\int_{\{(x,t)\in\R^{N+1}_+:t>\delta\}\setminus\mathcal
                                        K}
                                        \mathcal
                                        V_{\mathcal K,u}\mathcal V_{\mathcal K}
                                        \mathop{\rm div}(t^{1-2s}\nabla (u\eta\zeta_K))\,dx\,dt\bigg)\\
&\quad -\lim_{\delta\to0^+}\int_{\Omega\setminus K}\mathcal
                                        V_{\mathcal K,u}(x,\delta)u(x)\zeta_K(x) \delta^{1-2s} (\partial_t \mathcal V_{\mathcal K}) (x,\delta)\,dx\\
  &\quad +\lim_{\delta\to0^+}\int_{\partial\mathcal K\cap \{(x,t)\in\R^{N+1}_+:t>\delta\}}
    t^{1-2s} u\partial_\nu u\,dS \\
  &\quad +
    \lim_{\delta\to0^+}\int_{\partial\mathcal K\cap \{(x,t)\in\R^{N+1}_+:t>\delta\}}
                                       t^{1-2s} u^2\partial_\nu \mathcal V_{\mathcal K}\,dS.
\end{align*}
By Lemma \ref{l:prop_calV} and the Dominated Convergence Theorem, we
have
\begin{equation*}
  \lim_{\delta\to0^+}\int_{\Omega\setminus K}
  \mathcal V_{\mathcal K,u}(x,\delta)u(x)\zeta_K(x) \delta^{1-2s}
  (\partial_t \mathcal V_{\mathcal K}) (x,\delta)\,dx=0,
\end{equation*}
hence we conclude that
\begin{align}
  \label{eq:con1}
  \int_{\R^{N+1}_+\setminus\mathcal K}t^{1-2s}|\nabla \mathcal V_{\mathcal K,u}|^2\,dx\,dt&= 
                                        -2\int_{\R^{N+1}_+\setminus\mathcal K}
                                        t^{1-2s}\mathcal
                                        V_{\mathcal K,u}\nabla \mathcal V_{\mathcal K}
                                        \cdot\nabla
                                        (u\eta\zeta_K)\,dx\,dt\\
&\notag\quad -\int_{\R^{N+1}_+\setminus\mathcal K}
                                        \mathcal
                                        V_{\mathcal K,u}\mathcal V_{\mathcal K}
                                        \mathop{\rm div}(t^{1-2s}\nabla (u\eta\zeta_K))\,dx\,dt\\
 &\notag \quad+\int_{\partial\mathcal K\cap \R^{N+1}_+}
   t^{1-2s} u\partial_\nu u\,dS +
   \int_{\partial\mathcal K\cap \R^{N+1}_+}
                                       t^{1-2s} u^2\partial_\nu \mathcal V_{\mathcal K}\,dS.
\end{align}
Furthermore, an integration by parts yields
\begin{align}
\label{eq:con2} \int_{\mathcal K}t^{1-2s}|\nabla \mathcal V_{\mathcal K,u}|^2\,dx\,dt
  &= \int_{\mathcal K}t^{1-2s}|\nabla_xu|^2\,dx\,dt \\
 \notag&=
    -\int_{\mathcal K}t^{1-2s}u\Delta u\,dx\,dt-\int_{\partial\mathcal K\cap \R^{N+1}_+}
    t^{1-2s} u\partial_\nu u\,dS.
\end{align}
Combining \eqref{eq:con1} and \eqref{eq:con2}, we finally conclude
that
 \begin{align*}
    \mathcal C^s_\Omega(K,u)&=\frac{1}{\kappa_s}\int_{\R^{N+1}_+}t^{1-2s}|\nabla
   \mathcal V_{\mathcal K,u}|^2\,dx\,dt\\
   &=\frac{1}{\kappa_s}
\int_{\partial\mathcal K\cap\R^{N+1}_+}t^{1-2s}u^2\partial_\nu
    \mathcal V_{\mathcal K}\,dS- \frac{2}{\kappa_s} \int_{\R^{N+1}_+\setminus\mathcal
                                     K}t^{1-2s}\mathcal V_{\mathcal K,u}\nabla
    \mathcal V_{\mathcal K}\cdot\nabla (u\eta\zeta_K)\,dx\,dt\\
&\quad    -\frac{1}{\kappa_s}\int_{\R^{N+1}_+\setminus\mathcal K}\mathop{\rm
                                                                  div}(t^{1-2s}\nabla
                                                                  (u\eta\zeta_K))\mathcal
                                                                  V_{\mathcal
                                                                  K,u}
    \mathcal V_{\mathcal K}\,dx\,dt-\frac{1}{\kappa_s}\int_{\mathcal K}t^{1-2s}u\Delta u\,dx\,dt,
  \end{align*}
 thus completing the proof.
\end{proof}

By virtue of the approximation result in Lemma \ref{l:appr-cap_n},
passing to the limit in the identity established in Lemma
\ref{l:stimecap} yields the following proposition with a formula for
$ \mathop{\rm Cap}\nolimits^s_\Omega(K,u) $. In particular,
it is precisely from the boundary term
\begin{equation*}
\int_{\partial\mathcal K\cap\R^{N+1}_+}t^{1-2s}u^2\partial_\nu
\mathcal V_{\mathcal K}\,dS 
\end{equation*}
appearing in Lemma \ref{l:stimecap} that the $s$-equilibrium measure $\mu_{\Omega,K}^s$ introduced in Section
\ref{subsec_equil_measu} emerges upon taking the limit.

\begin{proposition}\label{p:for-cap}
  Let $\Omega\subset\R^N$ be a bounded open set and
  $s\in (0,\min\{1,N/2\})$.  Let $K\subset \Omega$ be a compact set,
  $\eta\in C^\infty_{\rm c}([0,+\infty))$ be such that $\eta(t)=1$ for
  all $t\in [0,1]$, and $\zeta_K\in C^\infty_c(\Omega)$ be such that
  $\zeta_K=1$ in a neighborhood of $K$. If $u\in C^\infty(\Omega)$,
  then
 \begin{align}
    \label{eq:lim-fo}    \mathop{\rm Cap}\nolimits^s_\Omega(K,u)=&
    \int_K u^2 \, d \mu_{\Omega,K}^s- \frac{2}{\kappa_s} \int_{\R^{N+1}_+}t^{1-2s}
                                            V_{\Omega,K,u}\nabla  V_{\Omega,K}\cdot\nabla (u\eta\zeta_K)\,dx\,dt\\
    \notag&   -\frac{1}{\kappa_s}\int_{\R^{N+1}_+}\mathop{\rm
            div}(t^{1-2s}\nabla
            (u\eta\zeta_K)) V_{\Omega,K,u}
            V_{\Omega,K}\,dx\,dt,
  \end{align}
 where $\mu_{\Omega,K}^s$ is the $s$-equilibrium measure of $K$
  introduced in Section \ref{subsec_equil_measu}.
\end{proposition}
\begin{proof}
  We first observe that there exists a sequence $\{\mathcal K_n\}_n$
  of compact sets in $\overline{\R^{N+1}_+}$ such that, for every $n$,
  $\mathcal K_n$ satisfies \eqref{eq:calK1}, \eqref{eq:calK3}, and
  \eqref{eq:calK2}, $\mathcal K_{n+1}\subseteq \mathcal K_{n}$, and
  $\bigcap_{n}\mathcal K_n=K \times\{0\}$.
 
Since $\zeta_K\in C^\infty_c(\Omega)$ is such that $\zeta_K=1$ in a
  neighborhood of $K$, we have $\zeta_K=1$ also in a
  neighborhood of $K_n:=\{(x,t)\in \mathcal K_n:t=0\}$ for $n$
  sufficiently large. Lemma \ref{l:stimecap} then applies to each
  $\mathcal K_n$, provided $n$ is sufficiently large,  and yields
  \begin{align}\label{eq:pass-limit}
    \mathcal C^s_\Omega(\mathcal K_n,u)=&
\frac{1}{\kappa_s}\int_{\partial\mathcal K_n\cap\R^{N+1}_+}t^{1-2s}u^2\partial_{\nu_n}  
                                                               \mathcal V_{\mathcal K_n}\,dS\\
   \notag &\quad- \frac{2}{\kappa_s} \int_{\R^{N+1}_+\setminus\mathcal
                                     K_n}t^{1-2s}\mathcal V_{\mathcal K_n,u}\nabla
    \mathcal V_{\mathcal K_n}\cdot\nabla (u\eta\zeta_K)\,dx\,dt\\
\notag&\quad    -\frac{1}{\kappa_s}\int_{\R^{N+1}_+\setminus\mathcal K_n}\mathop{\rm
                                                                  div}(t^{1-2s}\nabla
                                                                  (u\eta\zeta_K))\mathcal
                                                                  V_{\mathcal
                                                                  K_n,u}
        \mathcal V_{\mathcal K_n}\,dx\,dt\\
    \notag&\quad-\frac{1}{\kappa_s}\int_{\mathcal K_n}t^{1-2s}u\Delta u\,dx\,dt,
  \end{align}
  where $\nu_n$ is the unit vector on $\partial\mathcal K_n$ pointing
  toward the interior of $\mathcal K_n$.  
  By Lemma \ref{l:appr-cap_n}, we have
  \begin{equation*}
    \lim_{n\to\infty}  \mathcal C^s_\Omega(\mathcal K_n,u)
     =\mathop{\rm Cap}\nolimits^s_\Omega(K,u).
   \end{equation*}
   Moreover, again by  Lemma \ref{l:appr-cap_n}, $\mathcal V_{\mathcal
     K_n,u}\to V_{\Omega,K,u}$ and $\mathcal V_{\mathcal K_n}\to
   V_{\Omega,K}$  in $\mathcal{D}^{1,2}(\R^{N+1}_+;t^{1-2s})$. Then,
   taking into account that the $(N+1)$-dimensional Lebesgue measure
   of $\mathcal K_n$ tends to $0$ as $n\to\infty$, we have
   \begin{equation*}
     \lim_{n\to\infty} \int_{\R^{N+1}_+\setminus\mathcal
                                     K_n}t^{1-2s}\mathcal V_{\mathcal K_n,u}\nabla
                                   \mathcal V_{\mathcal K_n}\cdot\nabla (u\eta\zeta_K)\,dx\,dt
                                   =\int_{\R^{N+1}_+}t^{1-2s}V_{\Omega,K,u}\nabla
    V_{\Omega,K}\cdot\nabla (u\eta\zeta_K)\,dx\,dt
  \end{equation*}
  and
   \begin{equation*}
     \lim_{n\to\infty}\int_{\R^{N+1}_+\setminus\mathcal K_n}\mathop{\rm
                                                                  div}(t^{1-2s}\nabla
                                                                  (u\eta\zeta_K))\mathcal
                                                                  V_{\mathcal
                                                                  K_n,u}
    \mathcal V_{\mathcal K_n}\,dx\,dt=
\int_{\R^{N+1}_+}\mathop{\rm
                                                                  div}(t^{1-2s}\nabla
                                                                  (u\eta\zeta_K)) V_{\Omega,K,u}
    V_{\Omega,K}\,dx\,dt.
  \end{equation*}
  Furthermore, by absolute continuity of the integral,
  \begin{equation*}
    \lim_{n\to\infty}
    \int_{\mathcal K_n}t^{1-2s}u\Delta u\,dx\,dt=0.
  \end{equation*}
  Hence, passing  to the limit in  \eqref{eq:pass-limit} yields
  \begin{align*}
    \mathop{\rm Cap}\nolimits^s_\Omega(K,u)=
    &\mathcal L(u,K) - \frac{2}{\kappa_s} \int_{\R^{N+1}_+}t^{1-2s}
      V_{\Omega,K,u}\nabla  V_{\Omega,K}\cdot\nabla
      (u\eta\zeta_K)\,dx\,dt\\
    &   -\frac{1}{\kappa_s}\int_{\R^{N+1}_+}\mathop{\rm
      div}(t^{1-2s}\nabla
      (u\eta\zeta_K)) V_{\Omega,K,u}
      V_{\Omega,K}\,dx\,dt,
  \end{align*}
  with
  \begin{equation}\label{eq_limit_LuK}
    \mathcal L(u,K)=\frac{1}{\kappa_s}\lim_{n\to\infty}\int_{\partial\mathcal K_n\cap\R^{N+1}_+}t^{1-2s}u^2\partial_{\nu_n}  
    \mathcal V_{\mathcal K_n}\,dS.
  \end{equation}
 To complete the proof, it remains to show that
  \begin{equation}
    \label{eq_LmuK}
    \mathcal L(u,K)=\int_K u^2 \, d \mu_{\Omega,K}^s. 
  \end{equation}
  By \cite[Theorem 10.13]{L_book_sobolev}, for every $n\in\N$ there
  exists a unique Radon measure
  $\sigma_n: \mathcal B(\Omega) \to[0,+\infty]$ such that
\begin{equation}\label{eq_sigma_n}
  \int_\Omega \phi \, d \sigma_n=
  \frac1{\kappa_s}
  \int_{\partial  \mathcal K_n\cap\R^{N+1}_+}t^{1-2s}\phi\partial_{\nu_n}\mathcal V_{ \mathcal K_n}\,dS
\end{equation}
for any $\phi \in C^\infty_c(\Omega)$. We observe that
  \cite[Theorem 10.13]{L_book_sobolev} is applicable since the right
  hand side of \eqref{eq_sigma_n} defines a positive distribution in
view of Lemma \ref{l:prop_calV}.

For every compact set $F\subset\Omega$, let $\varphi_F\in
  C^\infty_c(\Omega)$ be such that $0\leq\varphi_F\leq1$ and
  $\varphi_F\equiv1$ on $F$. For every $n$, by \eqref{eq_sigma_n} and  Lemma
  \ref{l:prop_calV}-(vi) we have
  \begin{align*}
    \sigma_n(F)&\leq  \int_\Omega \varphi_F \, d \sigma_n=
     \frac1{\kappa_s}\int_{\partial  \mathcal
       K_n\cap\R^{N+1}_+}t^{1-2s}\varphi_F\partial_{\nu_n}\mathcal V_{
       \mathcal K_n}\,dS\\
     &\leq \frac1{\kappa_s} \int_{\partial  \mathcal
       K_n\cap\R^{N+1}_+}t^{1-2s}\partial_{\nu_n}\mathcal V_{ \mathcal
       K_n}\,dS=     \mathcal C^s_\Omega(\mathcal K_n),
  \end{align*}
  hence, by \eqref{eq:approx1}, $\sup_{n}\sigma_n(F)<+\infty$ for each
  compact set $F\subset\Omega$.  It follows that, up to subsequences,
  $\sigma_n \rightharpoonup \sigma$ weakly in the sense of measures as
  $n \to \infty$, for some positive Radon measure $\sigma$, see
  e.g. \cite[Theorem 1.41]{EG_book}. Furthermore, for any
  $\phi \in C^\infty_c(\Omega)$ and $\eta \in C^\infty_c([0,+\infty))$
  with $\eta(t)\equiv1$ in $[0,1]$, from an integration by parts,
  \eqref{eq:con-pot}, \eqref{eq:13}, and \eqref{eq:Lmu} it follows
  that
\begin{align*}
  \int_\Omega \phi \, d \sigma&=\lim_{n\to \infty}\int_\Omega \phi \, d \sigma_n
                                = \lim_{n\to \infty}\frac{1}{\kappa_s}\int_{\partial\mathcal K_n\cap\R^{N+1}_+}t^{1-2s}\phi \, \partial_{\nu_n}\mathcal V_{ \mathcal K_n}\,dS \\
                              &=\lim_{n\to \infty}\frac{1}{\kappa_s}\int_{\partial\mathcal K_n\cap\R^{N+1}_+}t^{1-2s}\phi\eta \, \partial_{\nu_n}\mathcal V_{ \mathcal K_n}\,dS
                                =\lim_{n\to \infty} \frac{1}{\kappa_s}q(\mathcal V_{\mathcal K_n},\phi
                                \eta)\\
                              &=\frac{1}{\kappa_s}q( V_{\Omega, K},\phi \eta)
                                =(v_{\Omega,K}, \phi)_{\mathcal D^{s,2}(\R^N)}
                                = \int_\Omega \phi \,d \mu_{\Omega,K}^s.
\end{align*}
We conclude that $\sigma= \mu_{\Omega,K}^s$ and, in fact, the
  full sequence $\sigma_n$ converges weakly to $\sigma$ in the sense
  of measures. By \eqref{eq_limit_LuK} and \eqref{eq_sigma_n} it
follows that
\begin{align*}
\mathcal L(u,K)&= \frac{1}{\kappa_s}\lim_{n\to\infty}\int_{\partial\mathcal K_n\cap\R^{N+1}_+}t^{1-2s}u^2\partial_{\nu_n}  
                 \mathcal V_{\mathcal K_n}\,dS\\
  &=\frac{1}{\kappa_s}
\lim_{n\to\infty}\int_{\partial\mathcal K_n\cap\R^{N+1}_+}t^{1-2s}\zeta_Ku^2\partial_{\nu_n}  
                 \mathcal V_{\mathcal K_n}\,dS\\
  &=\lim_{n\to\infty}\int_\Omega
    \zeta_Ku^2\,d\sigma_n=
    \int_\Omega \zeta_Ku^2 d \mu_{\Omega,K}^s=\int_K u^2 d \mu_{\Omega,K}^s,
\end{align*}
since $\mu_{\Omega,K}^s$ is supported on $K$. In conclusion, we have proved \eqref{eq_LmuK}.
\end{proof}

As a consequence of Proposition \ref{p:for-cap}, we obtain the
following asymptotic expansion for the fractional $u$-capacity along a
sequence of compact sets $\{K_n\}_{n\in\N}$ satisfying hypotheses
\eqref{eq:K1}, \eqref{eq:K2}, and \eqref{eq:K3}. The integral of $u^2$
with respect to the $s$-equilibrium measure of $K_n$ serves as  the leading
term of the expansion, while the other terms in \eqref{eq:lim-fo} turn
out to be negligible with respect to the $s$-capacity by virtue of
Lemma \ref{lem:trace_o_capacity}.

\begin{corollary}\label{cor:cor1}
  Let $\Omega\subset\R^N$ be a bounded open set and
  $s\in (0,\min\{1,N/2\})$.  Let $\{K_n\}_{n\in\N}$ be a sequence of
  compact sets contained in $\Omega$ satisfying \eqref{eq:K1},
  \eqref{eq:K2}, and \eqref{eq:K3}. If $u\in C^\infty(\Omega)$, then
  \begin{equation}\label{eq:exp-capa-n}
   \mathop{\rm Cap}\nolimits^s_\Omega(K_n,u)=
    \int_{K_n} u^2 \, d\mu_{\Omega,K_n}^s +o(\mathop{\rm
     Cap}\nolimits^s_\Omega(K_n))\quad\text{as }n\to\infty,
 \end{equation}
 where $\mu_{\Omega,K_n}^s$ is the $s$-equilibrium measure of $K_n$ as in
 Section \ref{subsec_equil_measu}. Furthermore, for every $n\in\N$,
  \begin{equation}\label{eq:cond_L}
\Big(\min_{K_n} u^2\Big) {\mathop{\rm
    Cap}\nolimits^s_\Omega(K_n)}\leq  \int_{K_n} u^2 \, d\mu_{\Omega,K_n}^s \leq \Big(\max_{K_n} u^2\Big) {\mathop{\rm
    Cap}\nolimits^s_\Omega(K_n)}.
  \end{equation}  
\end{corollary}
\begin{proof}
  By \eqref{eq_muK}, it is easy to check that \eqref{eq:cond_L} holds.

  Let us fix $\eta\in C^\infty_{\rm c}([0,+\infty))$ such that
  $\eta(t)=1$ for all $t\in [0,1]$.

  By assumption \eqref{eq:K3}, the set
  $K=\{x\in \Omega: \mathop{\rm dist}(x,\partial\Omega)\geq \delta\}\subset\Omega$
  is compact and $K_n\subseteq K$ for all $n\in\N$. Let
  $\zeta_K\in C^\infty_c(\Omega)$ be such that $\zeta_K\equiv1$ in a
  neighborhood of $K$. Then $\zeta_K\equiv1$ in a neighborhood of $K_n$ as
  well, for every $n\in\N$. Proposition
  \ref{p:for-cap} implies that, for every $n$,
  \begin{align}\label{eq:ptl1}
    \mathop{\rm Cap}\nolimits^s_\Omega(K_n,u)= \int_{K_n} u^2 \, d
    \mu_{\Omega,K_n}^s&
                      - \frac{2}{\kappa_s} \int_{\R^{N+1}_+}t^{1-2s}
                        V_{\Omega,K_n,u}\nabla  V_{\Omega,K_n}\cdot
                        \nabla (u\eta\zeta_K)\,dx\,dt\\
    \notag &   -\frac{1}{\kappa_s}\int_{\R^{N+1}_+}\mathop{\rm
             div}(t^{1-2s}\nabla
             (u\eta\zeta_K)) V_{\Omega,K_n,u}
             V_{\Omega,K_n}\,dx\,dt.
  \end{align}
  By assumption \eqref{eq:K2}, $\lim_{n\to\infty}\mathop{\rm
    Cap}\nolimits^s_\Omega(K_n)=0$. Hence Lemma \ref{lem:trace_o_capacity}
  applies: using \eqref{eq:rs} and \eqref{eq:o2} for some
  bounded open set $\mathcal U\subset\R^{N+1}_+$ whose closure
  contains the support of $(x,t)\mapsto \eta(t)\zeta_K(x)$, we obtain
  \begin{align*}
  &  \left|\int_{\R^{N+1}_+}t^{1-2s}
 V_{\Omega,K_n,u}\nabla  V_{\Omega,K_n}\cdot\nabla
    (u\eta\zeta_K)\,dx\,dt\right|\\
    &\quad\leq
    \Big(\sup|\nabla
    (u\eta\zeta_K)|\Big)
   \sqrt{ \int_{\mathcal U}t^{1-2s}
    V_{\Omega,K_n,u}^2\,dx\,dt}
\sqrt{\int_{\R^{N+1}_+}t^{1-2s}
      |\nabla  V_{\Omega,K_n}|^2\,dx\,dt}\\
    &\quad\leq
    \Big(\sup|\nabla
    (u\eta\zeta_K)|\Big)\Big(\max_{K_n}|u|\Big)
   \sqrt{ \int_{\mathcal U}t^{1-2s}
    V_{\Omega,K_n}^2\,dx\,dt}
\sqrt{\int_{\R^{N+1}_+}t^{1-2s}
      |\nabla  V_{\Omega,K_n}|^2\,dx\,dt}\\
  &\quad  =o(\mathop{\rm Cap}\nolimits^s_\Omega(K_n))
\quad\text{as } n\to\infty.
  \end{align*}
  In a similar way, we also obtain
  \begin{equation*}
\int_{\R^{N+1}_+}\mathop{\rm  div}(t^{1-2s}\nabla (u\eta\zeta_K)) V_{\Omega,K_n,u}
     V_{\Omega,K_n}\,dx\,dt  =o(\mathop{\rm
       Cap}\nolimits^s_\Omega(K_n))
     \quad\text{as } n\to\infty.
   \end{equation*}
   Expansion \eqref{eq:exp-capa-n} then directly follows from
   \eqref{eq:ptl1}.
\end{proof}
As observed in \cite[Example 2.5]{AFN}, the fractional capacity of a
point is zero, i.e. $\mathop{\rm Cap}\nolimits^s_\Omega(\{x_0\})=0$ for
every $x_0\in\Omega$.
Hence, if $\{K(\varepsilon)\}_{\varepsilon>0}$ is a family of compact
sets concentrating
at a singleton $\{x_0\}$, by \cite[Lemma 4.2]{AFN} we have that
$\lim_{\e\to 0}\mathop{\rm Cap}\nolimits^s_\Omega(K(\e))=0$.
Then, Corollary \ref{cor:cor1} applies, yielding the following result.
\begin{corollary}\label{cor:cor2}
  Let $\Omega\subset\R^N$ be a bounded open set and
  $s\in (0,\min\{1,N/2\})$.  For some $x_0\in\Omega$, let
  $\{K(\varepsilon)\}_{\varepsilon>0}$ be a family of compact sets
  contained in $\Omega$ such that $K(\varepsilon)$ is concentrating at
  $\{x_0\}$ in the sense of Definition \ref{def:concentrates}. If
  $u\in C^\infty(\Omega)$, then
  \begin{equation}\label{eq:expa5-7}
   \mathop{\rm Cap}\nolimits^s_\Omega(K(\e),u)=u^2(x_0) \mathop{\rm
     Cap}\nolimits^s_\Omega(K(\e))+o(\mathop{\rm
     Cap}\nolimits^s_\Omega(K(\e)))\quad\text{as }\e\to0.
 \end{equation}
\end{corollary}
\begin{proof}
  By the continuity of $u$ at $x_0$ and the concentration of the
  family $\{K(\e)\}_{\e>0}$ at $\{x_0\}$, we have that, for any
  sequence $\e_n\to0^+$, 
  \begin{equation*}
    \lim_{n\to\infty}\min_{K(\e_n)} u^2=
    \lim_{n\to\infty}\max_{K(\e_n)} u^2=u^2(x_0).
  \end{equation*}
  By \eqref{eq:cond_L} this implies that
  \begin{equation*}
    \int_{K(\e_n)} u^2 \, d\mu_{\Omega,K(\e_n)}^s=u^2(x_0)\mathop{\rm
     Cap}\nolimits^s_\Omega(K(\e_n))+o(\mathop{\rm
     Cap}\nolimits^s_\Omega(K(\e_n)))\quad\text{as }n\to\infty.
 \end{equation*}
Hence, from \eqref{eq:exp-capa-n} we deduce the asymptotic expansion
\eqref{eq:expa5-7} first along the sequence
 $\{\e_n\}$, and  then as $\e\to0$ by the Urysohn subsequence principle.
 \end{proof}

\section{Expansions of multiple eigenvalues}\label{sec_eigen}
In this section, we apply the \emph{Lemma on small eigenvalues} by
Y. Colin de Verdi\`ere \cite{ColindeV1986} to derive an asymptotic
expansion of the eigenvalue variation
$\lambda_{k}^s(\Omega\setminus K_n) -\lambda_k^s(\Omega)$, for
$\{K_n\}_{n\in\N}$ being a sequence of compact sets contained in
$\Omega$ satisfying \eqref{eq:K1}, \eqref{eq:K2}, \eqref{eq:K3}, and
\eqref{hp_weak_conv_mu}.  We focus on an eigenvalue
$\lambda_{n_0}^s(\Omega)$ with multiplicity $m_0\in \N\setminus\{0\}$,
i.e.  $n_0,m_0 \in \N\setminus\{0\}$ are such that
\eqref{lambda_mul_m} is satisfied.  Let $E(\lambda_{n_0}^s(\Omega))$
be the associated $m_0$-dimensional eigenspace, i.e.
\begin{equation*}
  E(\lambda_{n_0}^s(\Omega))=\big\{v\in \mathcal D^{s,2}(\Omega):
  v\text{ weakly solves
  }\eqref{eq:eigenvalue_main} \text{ with }\lambda=\lambda_{n_0}^s(\Omega)\big\}.
\end{equation*}
Let
\begin{equation}\label{eq:basis}
\{u_{n_0}, u_{n_0+1},\dots, u_{n_0+m_0-1}\} \quad\text{be a basis of $E(\lambda_{n_0}^s(\Omega))$ orthonormal in $L^2(\Omega)$}.
\end{equation}
For every $j\in\{1,\dots,m_0\}$, since
$u_{n_0+j-1} \in \mathcal{D}^{s,2}(\Omega)$ solves \eqref{eq_eigen}
(with extension $\mathcal H(u_{n_0+j-1})$ solving \eqref{eq:eige-exte}),
Schauder estimates for the fractional Laplacian ensure that
\begin{equation*}
u_{n_0+j-1} \in C^\infty(\Omega),
\end{equation*}
see e.g. \cite[Theorem 2.15]{TLX} and \cite{RS_reg_Lap}.

For any $u \in E(\lambda_{n_0}^s(\Omega))$ and $n\in\N$, let
$v_{\Omega,K_n}^u$ be as in \eqref{eq:pot-tr}. Let us define the
bilinear form
\begin{align}
\label{def_re}&r^n_{n_0}:E(\lambda_{n_0}^s(\Omega)) \times E(\lambda_{n_0}^s(\Omega)) \to  \R,\\
  &r^n_{n_0} (u,w):=\lambda_{n_0}^s(\Omega) \int_{\Omega}
    v_{\Omega,K_n}^u
    (w-v_{\Omega,K_n}^w) \, dx, \notag 
\end{align}
which will play an important role in the asymptotic expansion of $\lambda_{k}^s(\Omega\setminus K_n) -\lambda_k^s(\Omega)$.

\begin{remark}\label{remark_sym}
  The bilinear form $r^n_{n_0}$ is symmetric.
  Indeed, for any $u,w \in E(\lambda_{n_0}^s(\Omega))$ we have 
\begin{align*}
  r^n_{n_0} (u,w)-r^n_{n_0} (w,u)&=
                                   \lambda_{n_0}^s(\Omega)\int_{\Omega} wv_{\Omega,K_n}^u \, dx-
                                   \lambda_{n_0}^s(\Omega)\int_{\Omega} uv_{\Omega,K_n}^w\, dx\\
  &=(w,v_{\Omega,K_n}^u)_{\mathcal  D^{s,2}(\Omega)}-(u,v_{\Omega,K_n}^w)_{\mathcal D^{s,2}(\Omega)}\\
  &=(w-v_{\Omega,K_n}^w,v_{\Omega,K_n}^u)_{\mathcal  D^{s,2}(\Omega)}-(u-v_{\Omega,K_n}^u,v_{\Omega,K_n}^w)_{\mathcal D^{s,2}(\Omega)}\\
  &\quad +(v_{\Omega,K_n}^w,v_{\Omega,K_n}^u)_{\mathcal  D^{s,2}(\Omega)}-(v_{\Omega,K_n}^u,v_{\Omega,K_n}^w)_{\mathcal  D^{s,2}(\Omega)}=0,
\end{align*}
in view of \eqref{eq_eigen}  and \eqref{eq_vKu}. 
\end{remark}

\begin{proposition}\label{prop_re_capa}
Let $r^n_{n_0}$ be as in \eqref{def_re}. For any $u \in E(\lambda_{n_0}^s(\Omega))$, as $n\to\infty$
\begin{equation}\label{eq_re_capa}
r^n_{n_0} (u,u)=\mathop{\rm Cap}\nolimits^s_\Omega(K_n,u)+O(\|v^u_{\Omega,K_n}\|^2_{L^2(\Omega)}).
\end{equation}
\end{proposition}

\begin{proof}
By \eqref{def_re}, as $n\to\infty$,
\begin{equation*}
r^n_{n_0} (u,u)=\lambda_{n_0}^s(\Omega) \int_{\Omega} v_{\Omega,K_n}^u(u-v_{\Omega,K_n}^u) \, dx
=\lambda_{n_0}^s(\Omega) \int_{\Omega}u v_{\Omega,K_n}^u \, dx+O(\|v^u_{\Omega,K_n}\|^2_{L^2(\Omega)}).
\end{equation*}
Furthermore, by \eqref{eq_eigen}, \eqref{eq_vKu}, and \eqref{eq_capa_vKu} 
\begin{align*}
  \lambda_{n_0}^s(\Omega) \int_{\Omega}u v_{\Omega,K_n}^u \, dx
  &=(u,v_{\Omega,K_n}^u)_{\mathcal  D^{s,2}(\Omega)}\\
  &=(u-v_{\Omega,K_n}^u,v_{\Omega,K_n}^u)_{\mathcal  D^{s,2}(\Omega)}+\mathop{\rm Cap}\nolimits^s_\Omega(K_n,u)=\mathop{\rm Cap}\nolimits^s_\Omega(K_n,u).
\end{align*}
Hence, \eqref{eq_re_capa} is proved.
\end{proof}

  For every $n\in\N$, let
  \begin{equation*}
    \mathcal T_n:E(\lambda_{n_0}^s(\Omega))\to \mathcal D^{s,2}(\Omega\setminus K_n),\quad
    \mathcal T_n(u)=u-v_{\Omega,K_n}^u,
  \end{equation*}
  with $v_{\Omega,K_n}^u$ as in \eqref{eq:pot-tr}. 
\begin{proposition}\label{p:op-norm}
  For every $n\in\N$, let
  \begin{equation}\label{eq:def-tau-eps}
    \tau_n:=\|\mathop{\rm Id}-\mathcal T_n\|_{\mathcal
      L(E(\lambda_{n_0}^s(\Omega)),L^2(\Omega))},
  \end{equation}
  where $\mathop{\rm Id}$ is the identity map on $L^2(\Omega)$ and 
  $E(\lambda_{n_0}^s(\Omega))$ is meant as a subspace of $L^2(\Omega)$, i.e.
  \begin{equation}\label{def_taue}
   \tau_n=\sup\{\|u-\mathcal T_n u\|_{L^2(\Omega)}:
   u\in
   E(\lambda_{n_0}^s(\Omega)) \text{ and }\|u\|_{L^2(\Omega)}=1\}.
  \end{equation}
  Then
  \begin{equation}\label{eq:st-tau-eps}
    \tau_n\leq \sum_{j=1}^{m_0}\Big\|v_{\Omega,K_n}^{u_{n_0+j-1}}\Big\|_{L^2(\Omega)},
  \end{equation}
  where $\{u_{n_0+j-1}\}_{j=1}^{m_0}$ is a $L^2(\Omega)$-orthonormal basis
  of $E(\lambda_{n_0}^s(\Omega))$ as in \eqref{eq:basis}. Furthermore,
  \begin{equation}\label{eq_taue_smallo}
    \tau_n=o\big(\sqrt{\mathop{\rm Cap}\nolimits^s_\Omega(K_n)}\big)
  \end{equation}
  as $n\to\infty$.
\end{proposition}
\begin{proof}
  We have
  \begin{align*}
    \|\mathop{\rm Id}-\mathcal T_n\|_{\mathcal
    L(E(\lambda_{n_0}^s(\Omega)),L^2(\Omega))}&=
                                  \sup_{\substack{(\alpha_1,\dots,\alpha_{m_0})\in
                                  \R^{m_0}\\\sum_{j=1}^{m_0}
    \alpha_j^2=1}}\left\|(\mathop{\rm Id}-\mathcal T_n)\Big(\sum\nolimits_{j=1}^{m_0}
    \alpha_j u_{n_0+j-1}\Big)\right\|_{L^2(\Omega)}\\
    &=\sup_{\substack{(\alpha_1,\dots,\alpha_{m_0})\in
                                  \R^{m_0}\\\sum_{j=1}^{m_0}
    \alpha_j^2=1}} \left\|\sum\nolimits_{j=1}^{m_0}    \alpha_j
    v_{\Omega,K_n}^{u_{n_0+j-1}}\right\|_{L^2(\Omega)}\le
    \sum_{j=1}^{m_0}\Big\|v_{\Omega,K_n}^{u_{n_0+j-1}}\Big\|_{L^2(\Omega)},
  \end{align*}
  thus proving \eqref{eq:st-tau-eps}. Finally, by \eqref{eq:pot-tr},
  \eqref{eq:o1} and \eqref{eq:st-tau-eps} we obtain
  \eqref{eq_taue_smallo}.
\end{proof}

The following well-known \emph{Lemma on small eigenvalues}, due to Y. Colin
de Verdi\`ere \cite{ColindeV1986}, is an important tool in the proof
of our next result.

\begin{lemma}[Lemma on small eigenvalues \cite{ColindeV1986}]\label{lemma:CdV}
  Let $(H,(\cdot,\cdot)_{H})$ and
  $(\mathcal D,(\cdot,\cdot)_{\mathcal D})$ be real Hilbert spaces
  such that $\dim H=+\infty$, $\mathcal{D}$ is a dense subspace of
  $H$, and the embedding
  $\mathcal D \hookrightarrow \hookrightarrow H$ is continuous and
  compact.  Let $q: \mathcal{D}\times\mathcal{D}\to\R$ be a symmetric
  and  continuous bilinear form, such that $q$ is weakly coercive,
  i.e., there exist $A,B>0$ such that
  \begin{equation*}
    q(v,v)+A\|v\|_H^2\geq B\|v\|_{\mathcal D}^2\quad\text{for
      every }v\in \mathcal D,
  \end{equation*}
  where $\|\cdot\|_H=\sqrt{(\cdot,\cdot)_{H}}$ and
  $\|\cdot\|_{\mathcal D}=\sqrt{(\cdot,\cdot)_{\mathcal D}}$.

  Let $\{\nu_j\}_{j\in\N\setminus\{0\}}$ be the eigenvalues of $q$
  (repeated with their multiplicities, in ascending order) and let
  $\{{\mathbf e}_j\}_{j\in\N\setminus\{0\}}$ be an orthonormal basis
  of $H$ such that, for every $j\geq1$, ${\mathbf e}_j\in\mathcal D$
  is an eigenvector of $q$ associated to the eigenvalue $\nu_j$, i.e.
  $q({\mathbf e}_j,v)=\nu_j({\mathbf e}_j,v)_{H}$ for every
  $v\in \mathcal D$ and $j\geq1$.

  Let $m\in \N\setminus\{0\}$ and $F\subset \mathcal{D}$ be a
  $m$-dimensional linear subspace of $\mathcal D$. Let
  $\{\xi_j^F\}_{j=1,\dots,m}$ be the eigenvalues (in ascending order,
  repeated with their multiplicities) of $q$ restricted to $F$.

Assume that there exist $n\in\N\setminus\{0\}$ and
  $\gamma>0$ such that
	\begin{itemize}
        \item[\rm (H1)] $\nu_j\leq -\gamma$ for all $j\leq n-1$,
          $|\nu_j|\leq \gamma$ for all $j=n,\dots,n+m-1$, and 
          $\nu_j\geq \gamma$ for all $j\geq n+m$;
        \item[\rm (H2)]
          $\delta:=\sup\{|q(v,\varphi)|\colon v\in \mathcal{D},\ \varphi\in
          F,\ \|v\|_H=\|\varphi\|_H=1\}<\gamma/\sqrt{2}$.
	\end{itemize}
        Then
	\begin{equation*}
          \left|\nu_{n+j-1}-\xi_j^F\right|\leq\frac{4\delta^2}{\gamma}
          \quad\text{for all }j=1,\dots,m.
	\end{equation*}
\end{lemma}

We are going to apply Lemma \ref{lemma:CdV} in order to obtain an
asymptotic expansion for the eigenvalues variation.

\begin{proposition}\label{theo_esymp_ebstract}
  Let assumption \eqref{lambda_mul_m} be satisfied.  Let
  $\{K_n\}_{n\in\N}$ be a family of compact sets contained in $\Omega$
  satisfying \eqref{eq:K1}, \eqref{eq:K2}, \eqref{eq:K3}, and
  \eqref{hp_weak_conv_mu}.  For every $j\in\{1,\dots,m_0\}$, the
  $j$-th eigenbranch emanating from $\lambda_{n_0}^s(\Omega)$ can be
  expanded as
\begin{equation}\label{eq:expa1} 
  \lambda_{n_0+j-1}^s(\Omega\setminus K_n)=\lambda_{n_0}^s(\Omega)+\mu_j^n+O(\tau_n^2)\quad 
  \text{as }n\to\infty,
\end{equation}
where $\tau_n$ is defined in \eqref{eq:def-tau-eps} and
$\{\mu_j^n\}_{j=1}^{m_0}$ are the eigenvalues (taken in non-decreasing
order) of the bilinear form $r^n_{n_0}$ introduced in \eqref{def_re}.
\end{proposition}

\begin{proof}
We apply Lemma \ref{lemma:CdV} with the following choices:
\begin{align*}
  & \mathcal D=\mathcal D^{s,2}(\Omega\setminus K_n),\quad H=L^2(\Omega\setminus K_n),\\
  & q(u,v):=q_n(u,v)=(u,v)_{\mathcal D^{s,2}(\Omega\setminus K_n)}-\lambda^s_{n_0}(\Omega)\int_{\Omega}uv\,dx,\\
  & \nu_j=\lambda_j^s(\Omega\setminus K_n)-\lambda^s_{n_0} (\Omega),\quad F=\mathcal T_n(E(\lambda_{n_0}^s(\Omega))).
\end{align*}
We recall that the embedding
$\mathcal D^{s,2}(\Omega\setminus K_n) \hookrightarrow L^2(\Omega)$
is compact and, by Proposition \ref{p:op-norm}, we have that
$\|\mathop{\rm Id}-\mathcal T_n\|_{\mathcal
  L(E(\lambda_{n_0}^s(\Omega)),L^2(\Omega))}\to 0$ as $n\to\infty$.
Therefore, if $\{u_{n_0+j-1}\}_{j=1}^{m_0}$ is an $L^2(\Omega)$-orthonormal
basis of $E(\lambda_{n_0}^s(\Omega))$,
$\{\mathcal T_n(u_{n_0+j-1})\}_{j=1}^{m_0}$ are linearly independent in
$L^2(\Omega)$ as soon as $n$ is sufficiently large, so that
$F=\mathcal T_n(E(\lambda_{n_0}^s(\Omega)))$ is an $m_0$-dimensional
subspace of $\mathcal D^{s,2}(\Omega \setminus K_n)$.

Since $\lambda_{n_0}^s(\Omega)$ has multiplicity $m_0$ and, by  \eqref{limit_eigen}, 
$\lambda_k^s(\Omega\setminus K_n) \to \lambda_k^s(\Omega)$ as $n\to\infty$ for every $k\in\N\setminus\{0\}$, assumption (H1) in Lemma
\ref{lemma:CdV} holds whenever $n$ is large enough, if we set 
\begin{equation*}
  \gamma:={\frac12}\min\{\lambda_{n_0}^s(\Omega)-\lambda_{n_0-1}^s(\Omega),\lambda_{n_0+m_0}^s(\Omega)-\lambda_{n_0}^s(\Omega)\}
\end{equation*}
when $n_0\geq2$, while $\gamma:=\frac12
\left(\lambda_{2}^s(\Omega)-\lambda_{1}^s(\Omega)\right)$ when $n_0=1$.

Next, we bound the value $\delta$ in assumption (H2). For every
$v\in \mathcal D^{s,2}(\Omega\setminus K_n)$ and
$u\in E(\lambda_{n_0}^s(\Omega))$, by the Cauchy-Schwarz inequality,
\eqref{eq_eigen} and \eqref{eq_vKu} it follows that
\begin{align*}
|q_n(\mathcal T_n(u),v)|&= 
\left| (u-v_{\Omega,K_n}^u,v)_{\mathcal D^{s,2}(\Omega\setminus K_n)}-\lambda_{n_0}^s(\Omega)\int_{\Omega}  (u-v_{\Omega,K_n}^u)v\,dx\right|\\
&=\lambda_{n_0}^s(\Omega)\left|\int_\Omega  v_{\Omega,K_n}^u v\,dx\right|\\
&\leq\lambda_{n_0}^s(\Omega) \|(\mathop{\rm Id}-\mathcal
T_n)(u)\|_{L^2(\Omega)}\|v\|_{L^2(\Omega)}\\
& \leq\lambda_{n_0}^s(\Omega)\tau_n\|u\|_{L^2(\Omega)}\|v\|_{L^2(\Omega)}
\leq\lambda_{n_0}^s(\Omega)\frac{\tau_n}{1-\tau_n}\|\mathcal T_n (u)\|_{L^2(\Omega)}\|v\|_{L^2(\Omega)},
\end{align*}
where $\tau_n$ is defined in \eqref{eq:def-tau-eps}. The inequality
above is justified by the estimate
\begin{align*}
\|\mathcal T_n (u)\|_{L^2(\Omega)}&\ge\|u\|_{L^2(\Omega)}-\|(\mathop{\rm Id}-\mathcal T_n) (u)\|_{L^2(\Omega)}\\
&  \geq \|u\|_{L^2(\Omega)}-\tau_n\|u\|_{L^2(\Omega)}= (1-\tau_n)\|u\|_{L^2(\Omega)}.
\end{align*}
Therefore, 
\begin{align*}
  \delta:&=\sup\{|q_n(\mathcal T_n(u),v)|\colon v\in \mathcal D^{s,2}(\Omega\setminus K_n),u\in E(\lambda_{n_0}^s(\Omega)),\|v\|_{L^2(\Omega)}
           =\|\mathcal T_n(u)\|_{L^2(\Omega)}=1\}\\ 
\notag&\leq\lambda_{n_0}^s(\Omega) \frac{\tau_n}{1-\tau_n}
\end{align*}
and assumption (H2) of Lemma \ref{lemma:CdV} is satisfied if $n$ is
sufficiently large, in view of \eqref{eq_taue_smallo}.

Hence, we can apply Lemma \ref{lemma:CdV}, thus obtaining the
following expansion for every $j\in\{1,\dots,m_0\}$:
\begin{equation}\label{eq:asympt_proof}
\lambda_{n_0+j-1}^s(\Omega\setminus
K_n)=\lambda_{n_0}^s(\Omega)+\xi_j^n
+O(\tau_n^2)\quad\text{as }n\to\infty,
\end{equation}
where $\{ \xi_j^n\}_{j=1}^{m_0}$ are the eigenvalues (in
non-decreasing order) of the restriction of $q_n$ to
$\mathcal T_n(E(\lambda_{n_0}^s(\Omega)))$.

In order to characterize the values $\{\xi_j^n\}_{i=1,\dots,m_0}$
appearing in \eqref{eq:asympt_proof}, we consider the bilinear form
$r_{n_0}^n$ defined in \eqref{def_re} and observe that, by
\eqref{eq_eigen} and \eqref{eq_vKu},
\begin{align*}
  q_n(\mathcal T_n(u), \mathcal T_n(w))&=(u-v_{\Omega,K_n}^u,w-v_{\Omega,K_n}^w)_{\mathcal D^{s,2}(\Omega\setminus K_n)}-\lambda_{n_0}^s(\Omega)\int_\Omega (u-v_{\Omega,K_n}^u)(w-v_{\Omega,K_n}^w)\,dx\\
                                          &=(u,w-v_{\Omega,K_n}^w)_{\mathcal D^{s,2}(\Omega)}-(v_{\Omega,K_n}^u,w
                                            -v_{\Omega,K_n}^w)_{\mathcal D^{s,2}(\Omega)}\\&-  \lambda_{n_0}^s(\Omega)\int_\Omega u(w-v_{\Omega,K_n}^w)\,dx 
  +\lambda_{n_0}^s(\Omega)\int_\Omega v_{\Omega,K_n}^u(w-v_{\Omega,K_n}^w)\,dx\\
                                          &=(u,w-v_{\Omega,K_n}^w)_{\mathcal D^{s,2}(\Omega\setminus K_n)}\\
                                          &-\lambda_{n_0}^s(\Omega)\int_\Omega u(w-v_{\Omega,K_n}^w)\,dx+\lambda_{n_0}^s(\Omega)\int_\Omega v_{\Omega,K_n}^u(w-v_{\Omega,K_n}^w)\,dx\\
                                          &=\lambda_{n_0}^s(\Omega)\int_\Omega   v_{\Omega,K_n}^u(w-v_{\Omega,K_n}^w)\,dx=r^n_{n_0} (u,w)
\end{align*}
for every $u,w\in E(\lambda_{n_0}^s(\Omega))$.  Hence, the bilinear form
$q_n$ restricted to $\mathcal T_n(E(\lambda_{n_0}^s(\Omega)))$ coincides
with the bilinear form
\begin{equation*}
  (\phi,\psi)\in  \mathcal T_n(E(\lambda_{n_0}^s(\Omega)))
  \times \mathcal T_n(E(\lambda_{n_0}^s(\Omega))) 
  \mapsto r_{n_0}^n\big( \mathcal T_n^{-1}(\phi), \mathcal T_n^{-1} (\psi)\big).
\end{equation*}
Therefore, the eigenvalues $\{\xi_j^n\}_{j=1,\dots,m_0}$ of $q_n$
restricted to $\mathcal T_n(E(\lambda_{n_0}^s(\Omega)))$ coincide with
the eigenvalues of the $m_0\times m_0$ real symmetric matrix
\begin{equation}\label{eq:xi_mu_1}
\mathsf{B}_n=\mathsf{C}_n^{-1}\mathsf{R}_n,
\end{equation}
where
\begin{equation*}
\mathsf{R}_n=\big(r^n_{n_0} (u_{n_0+j-1}, u_{n_0+k-1})\big)_{jk}
\end{equation*}
and 
\begin{equation*}
\mathsf{C}_n=\big(c_{jk}^n \big)_{jk},\quad
c_{jk}^n=\int_\Omega \mathcal T_n(u_{n_0+j-1}) \mathcal T_n(u_{n_0+k-1})\,dx,
\end{equation*}
being   $\{u_{n_0+j-1}\}_{j=1}^{m_0}$ as in \eqref{eq:basis}. From
\eqref{eq:basis} and \eqref{eq:def-tau-eps} we have
\begin{equation*}
\mathsf{C}_n=\mathbb{I}+O(\tau_n)\quad\text{as }  n \to \infty,
\end{equation*}
with $\mathbb{I}$ being the identity $m_0\times m_0$ matrix. More
precisely, all the entries of the matrix $\mathsf{C}_n-\mathbb{I}$
are $O(\tau_n)$ as $n \to \infty$. Hence
$\mathsf{C}_n^{-1}=\mathbb{I}+O(\tau_n)$ as $n\to \infty$.  Moreover,
by \eqref{def_re} and \eqref{eq:def-tau-eps},
\begin{equation*}
\mathsf{R}_n=O(\tau_n)\quad\text{as }n \to \infty,
\end{equation*}
i.e.  all the entries of the matrix $\mathsf{R}_n$ are $O(\tau_n)$ as
$n \to \infty$.  As a consequence of \eqref{eq:xi_mu_1}, we obtain
\begin{equation}\label{eq:exp-matr}
  \mathsf{B}_n=(\mathbb{I}+O(\tau_n))\mathsf{R}_n=\mathsf{R}_n+O(\tau_n^2)
  \quad\text{as } n \to \infty.
\end{equation}
Since the  eigenvalues $\{\mu_j^n\}_{j=1}^{m_0}$ of the  bilinear form
$r^n_{n_0}$ are  the ones  of the  matrix $\mathsf{R}_n$,  the min-max
characterization  of eigenvalues  and \eqref{eq:exp-matr}  provide the
expansion
\begin{equation}\label{eq:expa-xi-mu}
\xi_j^n=\mu_j^n +O(\tau_n^2) \quad\text{as } n \to \infty,
\end{equation}
for all $j\in \{1,\dots,m_0\}$. Combining \eqref{eq:asympt_proof} and
\eqref{eq:expa-xi-mu} yields \eqref{eq:expa1}.
\end{proof}

We are now in a position to prove Theorem \ref{theo_main}.

\begin{proof}[Proof of Theorem \ref{theo_main}]
  Let us compute the asymptotics, as $n\to\infty$, of the bilinear
  form $r_{n_0}^n$ defined in \eqref{def_re}. By \eqref{eq_re_capa},
  \eqref{def_taue}, and \eqref{eq_taue_smallo}, for any
  $u \in E(\lambda_{n_0}^s(\Omega))$ we have
\begin{equation*}
  r^n_{n_0} (u,u)=\mathop{\rm Cap}\nolimits^s_\Omega(K_n,u)+O(\tau_n^2)=
  \mathop{\rm Cap}\nolimits^s_\Omega(K_n,u)
  +o(\mathop{\rm Cap}\nolimits^s_\Omega(K_n))
  \quad\text{as }n\to\infty.
\end{equation*}
By Corollary \ref{cor:cor1}   and \eqref{hp_weak_conv_mu}
\begin{align*}
  \mathop{\rm Cap}\nolimits^s_\Omega(K_n,u)
  &=\int_{K_n} u^2 \, d\mu_{\Omega,K_n}^s +o(\mathop{\rm Cap}\nolimits^s_\Omega(K_n))\\
  &=\left(\int_{\Omega} u^2 \, d\tilde\mu_{\Omega,K_n}^s\right)\mathop{\rm Cap}\nolimits^s_\Omega(K_n) +o(\mathop{\rm Cap}\nolimits^s_\Omega(K_n))\\
  &=\left(\int_{\Omega} u^2 \, d\tilde\mu\right)\mathop{\rm Cap}\nolimits^s_\Omega(K_n) +o(\mathop{\rm Cap}\nolimits^s_\Omega(K_n))
    \quad\text{as }n\to\infty.
\end{align*}
Hence, by the symmetry and bilinearity of $r^n_{n_0}$, for any
$u,v \in E(\lambda_{n_0}^s(\Omega))$,
\begin{equation*}
r^n_{n_0} (u,v)=\left(\int_{\Omega} u v \, d\tilde\mu\right)\mathop{\rm Cap}\nolimits^s_\Omega(K_n) +o(\mathop{\rm Cap}\nolimits^s_\Omega(K_n))
\quad\text{as }n\to\infty,
\end{equation*}
thus concluding the proof in view of \eqref{eq:expa1} and
\eqref{eq_taue_smallo}.
\end{proof}

We now apply Theorem  \ref{theo_main} to study splitting phenomena of multiple
eigenbranches, when perforating with holes that concentrate at a set of zero
capacity. In particular, we prove Theorem \ref{t:splitting}.

\begin{proof}[Proof of Theorem \ref{t:splitting}]
By contradiction, suppose that there exists a sequence $\e_n\to0^+$
such that 
\begin{equation}\label{eq:contr-ass}
  \lambda_{n_0}^{s}(\Omega\setminus K(\e_n))=
  \lambda_{n_0+1}^{s}(\Omega\setminus K(\e_n))=
  \dots=
  \lambda_{n_0+m_0-1}^{s}(\Omega\setminus K(\e_n))
\end{equation}
for all $n\in\N$.  By \cite[Lemma 4.2]{AFN}, we have that
$\lim_{n\to \infty}\mathop{\rm Cap}\nolimits^s_\Omega(K(\e_n))=0$.
Moreover, there exists some compact set $\bar K\subset\Omega$ such that
$K(\e_n)\subset\bar K$ for $n$ sufficiently large. On the other hand,
by \cite[Theorem 1.41]{EG_book}, along some subsequence still denoted
by $\{\e_n\}$, there exists some probability Radon measure $\tilde\mu$
with support contained in $K$ such that
$\tilde{\mu}^s_{\Omega,K(\e_{n})}$ weakly converges to $\tilde \mu$.
Hence, the sequence of compact sets $\{K(\e_n)\}$ satisfies conditions
\eqref{eq:K1}, \eqref{eq:K2}, \eqref{eq:K3}, and
\eqref{hp_weak_conv_mu}, so that Theorem \ref{theo_main} yields the
expansion
\begin{equation}\label{eq:expa-th-24}
  \lambda_{n_0+i-1}^{s}(\Omega\setminus K(\e_n)) =\lambda^s_{n_0} (\Omega)+
  \beta_i\mathop{\rm Cap}\nolimits^s_\Omega(K(\e_n))+
  o\left(\mathop{\rm Cap}\nolimits^s_\Omega(K(\e_n))\right)
  \quad\text{as } n\to\infty,
\end{equation}
for every $1\leq i\leq m_0$, with $\{\beta_i\}_{i=1,\dots,m_0}$ being the
eigenvalues (in increasing order) of the bilinear form $r_{\tilde\mu}$
defined in \eqref{eq:bilform}, which coincide with the eigenvalues of
the $m_0\times m_0$ matrix
\begin{equation*}
\mathsf{A}=\big(a_{jk} \big)_{jk},\quad
a_{jk}=\int_\Omega u_{n_0+j-1}u_{n_0+k-1} \, d\tilde \mu
=\int_K u_{n_0+j-1}u_{n_0+k-1} \, d\tilde \mu.
\end{equation*}
Assumption \eqref{eq:ass-conn-comp} implies that
\begin{equation*}
  a_{jk}\neq0\quad\text{for all }j,k\in\{1,\dots,m_0\}.
\end{equation*}
Since the eigenvalues of a real symmetric square matrix are all equal
to each other if and only if the matrix is a multiple of the identity
matrix, we conclude that 
 there exist
$j,k \in \{1,2,\dots,m_0\}$ such that $\beta_{j}\neq \beta_k$. In view
of \eqref{eq:expa-th-24} this implies that
$\lambda_{n_0+j-1}^{s}(\Omega\setminus K(\e_n))
\neq\lambda_{n_0+k-1}^{s}(\Omega\setminus K(\e_n))$ provided $n$ is
sufficiently large, thus contradicting  \eqref{eq:contr-ass}.
\end{proof}

We conclude by proving the fractional counterpart of Flucher's classic
result \cite[Theorem 6]{flucher}, namely Theorem \ref{cor_main},
dealing with a family of compact sets concentrating at a singleton
$\{x_0\}$.

\begin{proof}[Proof of Theorem \ref{cor_main}]
  Let $x_0 \in \Omega$ and $\{K(\varepsilon)\}_{\varepsilon>0}$ be a
  family of compact sets in $\Omega$ concentrating at $\{x_0\}$.
  Since $\mathop{\rm Cap}\nolimits^s_\Omega(\{x_0\})=0$, see
  \cite[Example 2.5]{AFN}, by \cite[Lemma 4.2]{AFN} we have
  $\lim_{\e\to 0}\mathop{\rm Cap}\nolimits^s_\Omega(K(\e))=0$.
  Moreover, we can assume without loss of generality that
  $\mathop{\rm Cap}\nolimits^s_\Omega(K(\e))>0$ for $\e>0$
  sufficiently small (otherwise, the expansions
  \eqref{eq:cor_main1}-\eqref{eq:cor_main2} are obviously
  satisfied). Since $x_0$ lies in the interior of $\Omega$, there
  exist some compact set $\bar K\subset\Omega$ such that
  $K(\e)\subset\bar K$ for $\e$ sufficiently small.

  By \cite[Theorem 1.41]{EG_book}, for any sequence $\e_n\to0^+$,
  there exists a subsequence $\{\e_{n_k}\}_{k\in\N}$ such that
  $\tilde{\mu}^s_{\Omega,K(\e_{n_k})}$ weakly converges to some
  probability Radon measure $\tilde\mu$ with support $\{x_0\}$, so
  that, necessarily, $\tilde\mu=\delta_{x_0}$. In particular, the
  sequence $\{K(\e_{n_k})\}_{k\in\N}$ satisfies assumptions \eqref{eq:K1},
  \eqref{eq:K2}, \eqref{eq:K3}, and \eqref{hp_weak_conv_mu}.  Letting
  $\{u_{n_0+j-1}\}_{j=1}^{m_0}$ be as in \eqref{eq:basis}, Theorem
  \ref{theo_main} then yields the following expansions:
\begin{equation*}
  \lambda_{n_0+i-1}^{s}(\Omega\setminus K(\e_{n_k})) =\lambda^s_{n_0} (\Omega)+
  \beta_i\mathop{\rm Cap}\nolimits^s_\Omega(K(\e_{n_k}))+
  o\left(\mathop{\rm Cap}\nolimits^s_\Omega(K(\e_{n_k}))\right)\quad\text{as } k\to\infty,
\end{equation*}
for every $i=1,\dots,m_0$, where $\{\beta_i\}_{i=1,\dots,m_0}$ are the
eigenvalues (in increasing order) of the bilinear form
  \begin{equation*}
(u,v)\in E(\lambda^s_{n_0}(\Omega))\times E(\lambda^s_{n_0}(\Omega)) \
\mapsto\ 
\int_\Omega u v \, d\delta_{x_0}=u(x_0)v(x_0),
\end{equation*}
i.e. $\{\beta_i\}_{i=1,\dots,m_0}$ are the eigenvalues of the
$m_0 \times m_0$-matrix
\begin{equation*}
\mathsf{A}=\big(u_{n_0+j-1}(x_0)u_{n_0+\ell-1}(x_0)\big)_{j\ell}.
\end{equation*}
It is simple linear algebra to see that the eigenvalues of
$\mathsf{A}$ are $0$ with multiplicity $m_0-1$ and
$\sum_{i=1}^{m_0}u_{n_0+i-1}^2(x_0)$ with multiplicity $1$, so that
\begin{equation*}
  \beta_1=\cdots=\beta_{m_0-1}=0,\quad\beta_{m_0}
  =\sum_{i=1}^{m_0}u_{n_0+i-1}^2(x_0).
\end{equation*}
Hence, the expansions \eqref{eq:cor_main1}-\eqref{eq:cor_main2} are
proved along the subsequence $\{\e_{n_k}\}_{k\in\N}$. Since the
coefficients $\{\beta_i\}_{i=1}^{m_0}$ depend neither on the sequence
$\{\e_n\}_{n\in\N}$ nor on the subsequence $\{\e_{n_k}\}_{k\in\N}$, by
the Urysohn subsequence principle we conclude that
\eqref{eq:cor_main1} and \eqref{eq:cor_main2} hold true as $\e\to0$.
\end{proof}

\section*{Acknowledgements}
\noindent
V.~Felli and G.~Siclari are partially supported by the 2026
INdAM--GNAMPA project ``Asymptotic analysis of variational problems''
(CUP E53C25002010001). P.~Roychowdhury is supported by the Seed Grant
provided by the Indian Institute of Technology Hyderabad under Grant
No.~SG/IITH/F391/2026-27/SG-230.

\bibliographystyle{acm}
\bibliography{bibliography}	

\end{document}